\documentclass[11pt]{article}
\usepackage[english]{babel}
\usepackage{lineno} 
\usepackage{graphicx,multicol}
\usepackage{epic,eepic,epsfig}
\usepackage{caption}

\usepackage{amsmath,amsfonts,amsthm,amssymb}
\usepackage{pstricks,pstricks-add,pst-math,pst-xkey}

\usepackage{algorithm, algorithmic}

\usepackage{graphicx,pstricks,pst-node,amssymb}

\usepackage{tikz}
\pgfdeclarelayer{edgelayer}
\pgfdeclarelayer{nodelayer}
\pgfsetlayers{edgelayer,nodelayer,main}

\tikzstyle{none}=[inner sep=0pt]
\definecolor{hexcolor0xf81e1c}{rgb}{0.973,0.118,0.110}
\definecolor{hexcolor0x3c00ff}{rgb}{0.235,0.000,1.000}

\tikzstyle{whitevertex}=[circle,fill=white,draw=black, scale = 0.5]
\tikzstyle{redvertex}=[circle,fill=hexcolor0xf81e1c,draw=black, scale = 0.5]
\tikzstyle{bluevertex}=[circle,fill=hexcolor0x3c00ff,draw=black, scale = 0.5]
\tikzstyle{greenvertex}=[circle,fill=green,draw=black, scale=0.5]
\tikzstyle{purplevertex}=[circle,fill=magenta,draw=black, scale=0.5]
\tikzstyle{grayvertex}=[circle,fill=white,draw=gray, scale=0.5]
\tikzstyle{blackvertex}=[circle,fill=black,draw=black, scale=0.5]

\tikzstyle{textbox}=[rectangle,fill=none,draw=none]
\tikzstyle{box}=[rectangle,fill=none,draw=black]

\tikzstyle{arc}=[black, ->]
\tikzstyle{grayarc}=[gray, ->]
\tikzstyle{bluearc}=[blue, ->]
\tikzstyle{grayedge}=[draw=gray]
\tikzstyle{blueedge}=[draw=blue]
\tikzstyle{rededge}=[draw=red]
\tikzstyle{edge}=[draw=black]

\tikzstyle{vertex}=[circle, ,fill=white,draw=black, scale=0.5]
\tikzstyle{10circle}=[circle, scale=10.0,draw=black]
\tikzstyle{10oval}=[ellipse, scale=10.0,draw=black]

\usepackage{amsmath,amsfonts,amsthm,amssymb}
\usepackage{pstricks,pstricks-add,pst-math,pst-xkey}

\usepackage{algorithm, algorithmic}

\usepackage{graphicx,pstricks,pst-node,amssymb}

\usepackage{algorithm, algorithmic}

\begin{document}

\newtheorem{theorem}{\hspace{5mm}Theorem}[section]
\newtheorem{prp}[theorem]{\hspace{5mm}Proposition}
\newtheorem{definition}[theorem]{\hspace{5mm}Definition}
\newtheorem{lemma}[theorem]{\hspace{5mm}Lemma}
\newtheorem{cor}[theorem]{\hspace{5mm}Corollary}
\newtheorem{alg}[theorem]{\hspace{5mm}Algorithm}
\newtheorem{sub}[theorem]{\hspace{5mm}Algorithm}
\newcommand{\induce}[2]{\mbox{$ #1 \langle #2 \rangle$}}
\newcommand{\2}{\vspace{2mm}}
\newcommand{\dom}{\mbox{$\rightarrow$}}
\newcommand{\ndom}{\mbox{$\not\rightarrow$}}
\newcommand{\compdom}{\mbox{$\Rightarrow$}}
\newcommand{\cdom}{\compdom}
\newcommand{\sdom}{\mbox{$\Rightarrow$}}
\newcommand{\lsd}{locally semicomplete digraph}
\newcommand{\lt}{local tournament}
\newcommand{\la}{\langle}
\newcommand{\ra}{\rangle}
\newcommand{\pf}{{\bf Proof: }}
\newtheorem{claim}{Claim}
\newcommand{\beq}{\begin{equation}}
\newcommand{\eeq}{\end{equation}}
\newcommand{\<}[1]{\left\langle{#1}\right\rangle}

\newcommand{\Z}{\mathbb{$Z$}}
\newcommand{\Q}{\mathbb{$Q$}}
\newcommand{\R}{\mathbb{$R$}}

%\bibliographystyle{plain}
%\pagewiselinenumbering
%\setpagewiselinenumbers
%\modulolinenumbers[1]
%\linenumbers

\title{Chordal signed graphs and signed bigraphs}

\author{Jing Huang\thanks{Department of Mathematics and Statistics, 
University of Victoria, Victoria, B.C., Canada; huangj@uvic.ca
(Research supported by NSERC)}\ \ and\ 
        Ying Ying Ye\thanks{Department of Mathematics and Statistics, 
University of Victoria, Victoria, B.C., Canada; fayye@uvic.ca}}

\date{}

\maketitle

\begin{abstract}
Chordal graphs and chordal bigraphs enjoy beautiful characterizations, in terms of 
forbidden subgraphs, vertex/edge orderings, vertex/edge separating sets, and 
tree-like representations. 
In this paper, we introduce chordal signed graphs and chordal signed bigraphs. 
Interestingly, chordal signed graphs are equivalent to strict chordal digraphs 
studied by Hell and Hern\'andez-Cruz. A forbidden subdigraph characterization of 
strict chordal digraphs can be translated to a forbidden subgraph characterization of
chordal signed graphs. We give a forbidden subgraph characterization of chordal 
signed bigraphs. The forbidden subgraphs for chordal signed bigraphs are analogous to 
those for chordal signed graphs but the proofs are much more complicated and 
intriguing.
\end{abstract}

\section{Introduction}

A graph is {\em chordal} if it does not contain an induced cycle of length $\geq 4$. 
Chordal graphs can be equivalently defined in terms of vertex orderings.
A vertex $v$ in a graph $G$ is {\em simplicial} if $N(v)$ induces a clique in $G$. 
A {\em perfect elimination ordering} of $G$ is a vertex ordering 
$v_1, v_2, \dots, v_n$ such that each $v_i$ is a simplicial vertex in 
$G-\{v_1, v_2, \dots, v_{i-1}\}$. A graph is chordal if and only if it has a 
perfect elimination ordering \cite{dirac,golumbic} 

A bipartite graph is {\em chordal} if it does not contain an induced cycles of length
$\geq 6$. Similar to chordal graphs, chordal bipartite graphs can be equivalently
defined in terms of edge orderings. 
Let $G$ be a bipartite graph with bipartition $(X,Y)$. A subgraph $H$ of $G$ is 
called a {\em biclique} of $G$ if every vertex in $V(H)\cap X$ is adjacent to all 
vertices in $V(H)\cap Y$. An edge $uv$ in $G$ is {\em simplicial} if $N(u) \cup N(v)$
induces a biclique in $G$. A bipartiate graph is chordal if and only if its edges
can be ordered $e_1, e_2, \dots, e_m$ in such a way that each $e_i$ is a simplicial 
edge in $G-\{e_1, e_2, \dots, e_{i-1}\}$ \cite{golumbic,gg}. Such an edge ordering 
is called as a {\em perfect edge-without-vertex elimination ordering} of $G$ 
\cite{bls}.   

Chordal graphs and chordal bipartite graphs are important in structural graph theory.
They admit elegant characterizations, in terms of, in addition to forbidden subgraphs
and perfect elimination orderings, vertex/edge separating sets and tree-like 
representations, which can be used to solve many optimization problems efficiently
\cite{golumbic,gg,hoc,huang,kk,klp,rtl}.

A {\em signed graph} $\widehat{G}$ is a graph $G$ in which each edge is signified
with either a positive sign or a negative sign \cite{harary}. 
Signed graphs have applications in geometry, matroid theory, and social
science \cite{harary,haraka,zas1,zas,zasla,zaslav}.
In a more general context an edge in a signed graph may have both positive and 
negative signs but this is not the case in this paper.

Let $\widehat{G}$ be a signed graph. An edge in $\widehat{G}$ is called {\em positive}
if it has a positive sign and {\em negative} otherwise.
A {\em subgraph} of $\widehat{G}$ is a signed graph $\widehat{H}$ obtained from 
$\widehat{G}$ by deleting vertices and/or edges. If $\widehat{H}$ is obtained by 
vertex deletions only then it is called an {\em induced subgraph} of $\widehat{G}$
or a subgraph of $\widehat{G}$ {\em induced by} the vertex set $V(\widehat{H})$ of 
$\widehat{H}$. We call a subgraph of $\widehat{G}$ {\em positive} if its edges are 
all positive.

We extend the concepts of chordal graphs and chordal bipartite graphs to signed
graphs.

We call a vertex $v$ in a signed graph $\widehat{G}$ {\em signed simplicial} if $N(v)$
induces a positive clique in $\widehat{G}$. We call a signed graph $\widehat{G}$ 
{\em chordal} if every induced subgraph of $\widehat{G}$ has a signed simplicial 
vertex, or equivalently, the vertices of $\widehat{G}$ can be ordered 
$v_1, v_2, \dots, v_n$ such that each $v_i$ is signed simplicial in 
$\widehat{G}-\{v_1,v_2,\dots,v_{i-1}\}$. 

Interestingly, an equivalent class of chordal signed graphs has been studied by
Hell and Hern\'andez-Cruz \cite{2h} under the name of 
{\em strict chordal digraphs}\footnote{Strict chordal digraphs are a subclass of 
chordal digraphs studied in \cite{hr,hy,meister}}.
Let $D$ be a digraph and let $G$ be the underlying graph of $D$. Form a signed graph 
$\widehat{G}$ by signifying edges $uv$ of $G$ with positive sign if $u$ and $v$ are 
joined with symmetric arcs and with negative sign if $u$ and $v$ are joined with 
a non-symmetric arc. Strict chordal digraphs studied in \cite{2h} are exactly 
the digraphs for which the associated signed graphs are chordal.
Hell and Hern\'andez-Cruz \cite{2h} gave a forbidden subdigraph characterization of 
strict chordal digraphs. This characterization yields a forbidden subgraph 
characterization of chordal signed graphs via the above transformation from digraphs 
to signed graphs.

For simplicity we shall use {\em bigraphs} for bipartite graphs. Thus,
a {\em signed bigraph} is a signed graph $\widehat{G}$ where $G$ is bipartite. 
When $G$ is a chordal bigraph, $\widehat{G}$ is called a {\em signed chordal bigraph}.
 
Let $\widehat{G}$ be a signed bigraph and $uv$ be an edge of $\widehat{G}$. 
Denote $N(uv) = (N(u) \cup N(v)) - \{u,v\}$. We call $uv$ {\em signed simplicial} if
$N(uv)$ induces a positive biclique in $\widehat{G}$. Note that a signed simplicial 
edge may be adjacent to negative edges.
A {\em signed simplicial edge-without-vertex ordering} of $\widehat{G}$ is an edge 
ordering $e_1, e_2, \dots, e_m$ of $\widehat{G}$ such that each $e_i$ is 
signed simplicial in $\widehat{G}-\{e_1,e_2,\dots,e_{i-1}\}$. 
We call $\widehat{G}$ a {\em chordal signed bigraph} if it has a signed simplicial 
edge-without-vertex ordering.

It follows from definition that a signed simplicial edge in $\widehat{G}$ is 
a simplicial edge in $G$.
Thus if $\widehat{G}$ is a chordal signed bigraph then $G$ is a chordal bigraph. 
Not all signed chordal graphs are chordal signed graphs. For example, the negative
$C_4$ is a signed chordal graph but not a chordal signed graph. On the other hand,  
any positive chordal bigraph is a chordal signed bigraph.
It also follows from definition that an edge in $\widehat{G}$ with one of its 
endvertices having degree one is a signed simplicial edge. Consequently,
any signed tree is a chordal signed bigraph. 

The main result of this paper is a forbidden subgraph characterization of chordal 
signed bigraphs (see Theorem \ref{main}). A graph is {\em non-trivial} if it has 
at least one edge. Our characterization implies that a signed graph is not chordal 
if and only if it contains a non-trivial induced subgraph which has no signed 
simplicial edge (see Corollary \ref{maincor}).

We will use edge-coloured graphs to represent signed graphs where positive edges are
coloured {\em blue}, negative edges are coloured {\em red}, and edges whose signs 
are not specified are coloured {\em black}. 
Suppose that $F$ is an edge-coloured graph whose edges are coloured with blue, red, 
and black colours. If $F$ has no black edge then it represents a unique signed 
graph. If $F$ has black edges then by replacing each black edge of $F$ with a blue 
or a red edge we obtain a signed graph and thus $F$ represents a set of signed 
graphs. 

The paper is organized as follows: In Section \ref{2}, we focus on complete bigraphs
and characterize chordal signed complete bigraphs by forbidden subgraphs 
(see Theorem \ref{comp}). 
In Section \ref{3}, we extend the characterization of chordal signed complete bigraphs
to a forbidden subgraph characterization of chordal signed non-separable bigraphs
(see Theorem \ref{non-separableForbidden}).  
In Section \ref{4}, we analyze the structure of chordal signed separable bigraphs and
find their forbidden subgraphs. Finally, in Section \ref{5}, we prove the main 
theorem (that is, Theorem \ref{main}).

\section{Complete bigraphs} \label{2}

In this section we consider signed complete bigraphs. We will show that there are
only five minimal signed complete bigraphs which are not chordal 
(see Figure~\ref{minForb-complete}).

\begin{lemma}\label{comp-lem}
Let $\widehat{G}$ be a signed complete bigraph with bipartition $(X,Y)$. 
If $\widehat{G}$ has a signed simplicial edge then it is a chordal signed bigraph.
\end{lemma}
\pf Denote $X = \{x_1, x_2, \dots, x_\alpha\}$ and $Y = \{y_1, y_2, \dots, y_\beta\}$.
Assume without loss of generality that $x_1y_1$ is a signed simplicial edge. 
Then $\widehat{G} - \{x_1,y_1\}$ is a positive biclique. It follows that
\[x_1y_1, x_1y_2, \dots, x_1y_\beta, x_2y_1, x_2y_2, \dots, x_2y_\beta, \dots, 
x_{\alpha}y_1, x_{\alpha}y_2, \dots, x_{\alpha}y_{\beta} \]
is a signed simplcial edge-without-vertex ordering of $\widehat{G}$. Hence
$\widehat{G}$ is a chordal singed bigraph.
\qed

\begin{center}
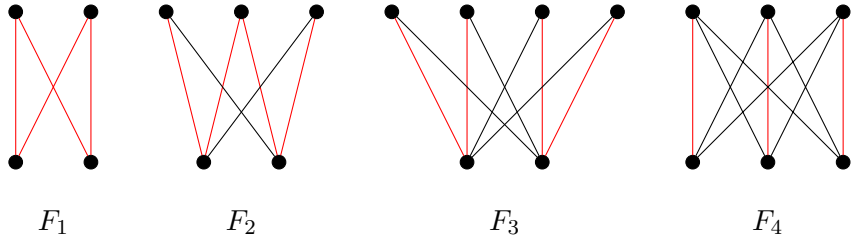
\begin{figure}[htb]
                \center
                \begin{tikzpicture}[>=latex]
                \begin{pgfonlayer}{nodelayer}
                %Fig a
\node [style=blackvertex] (1) at (0,0) {};
\node [style=blackvertex] (2) at (1,0) {};
\node [style=blackvertex] (3) at (0,2) {};
\node [style=blackvertex] (4) at (1,2) {};

\draw[style=edge] (1) to (3) to (2) to (4) to (1) [red];
\node [style=textbox] at (.5,-.8) {$F_1$};

\node [style=blackvertex] (5) at (2.5,0) {};
\node [style=blackvertex] (6) at (3.5,0) {};
\node [style=blackvertex] (7) at (2,2) {};
\node [style=blackvertex] (8) at (3,2) {};
\node [style=blackvertex] (9) at (4,2) {};

\draw[style=edge] (7) to (5) to (8) to (6) to (9)[red];
\draw[style=edge] (6) to (7);
\draw[style=edge] (5) to (9);
\node [style=textbox] at (3,-.8) {$F_2$};

\node [style=blackvertex] (10) at (6,0) {};
\node [style=blackvertex] (11) at (7,0) {};
\node [style=blackvertex] (12) at (5,2) {};
\node [style=blackvertex] (13) at (6,2) {};
\node [style=blackvertex] (14) at (7,2) {};
\node [style=blackvertex] (15) at (8,2) {};

\draw[style=edge] (12) to (10) to (13)[red];
\draw[style=edge] (14) to (11) to (15)[red];
\draw[style=edge] (12) to (11) to (13);
\draw[style=edge] (14) to (10) to (15);
\node [style=textbox] at (6.5,-.8) {$F_3$};

\node [style=blackvertex] (16) at (9,0) {};
\node [style=blackvertex] (17) at (10,0) {};
\node [style=blackvertex] (18) at (11,0) {};
\node [style=blackvertex] (19) at (9,2) {};
\node [style=blackvertex] (20) at (10,2) {};
\node [style=blackvertex] (21) at (11,2) {};

\draw[style=edge] (16) to (19)[red];
\draw[style=edge] (17) to (20)[red];
\draw[style=edge] (18) to (21)[red];
\draw[style=edge] (16) to (20) to (18) to (19) to (17) to (21) to (16);
\node [style=textbox] at (10,-.8) {$F_4$};

                \end{pgfonlayer}
                \end{tikzpicture}
\caption{\label{Forb-complete}The forbidden subgraphs for chordal signed complete 
bigraphs}
        \end{figure}
\end{center}

Figure~\ref{Forb-complete} depicts four sets $F_1, F_2, F_3, F_4$ of signed graphs.
It is easy to verify that no graph in $F_1 \cup F_2 \cup F_3 \cup F_4$ has a signed 
simplicial edge and hence none of them can be an induced subgraph of a chordal signed 
bigraph.

\begin{theorem} \label{comp}
Let $\widehat{G}$ be a signed complete bigraph. Then $\widehat{G}$ is chordal 
if and only if it does not contain any graph in $F_1 \cup F_2 \cup F_3 \cup F_4$ as 
an induced subgraph.
\end{theorem}
\pf The necessity of the statement follows from the remark above. 
For the sufficiency suppose to the contrary that there exists a signed complete 
bigraph which does not contain any graph in $F_1 \cup F_2 \cup F_3 \cup F_4$ as 
an induced subgraph and is not chordal. Let $\widehat{G}$ be a minimal such bigraph. 
Since $\widehat{G}$ is not chordal, it is non-trivial and by Lemma \ref{comp-lem} it 
has no signed simplicial edge.

Let $(X,Y)$ be a bipartition of $\widehat{G}$ where 
$X = \{x_1, x_2, \dots, x_{\alpha}\}$ and $Y = \{y_1, y_2, \dots, y_{\beta}\}$.
We show first that $\widehat{G}$ contains no negative $P_4$. Suppose to the contrary 
that $x_1y_1x_2y_2$ is a negative $P_4$ in $\widehat{G}$. Since no edge in 
$\widehat{G}$ is signed simplicial, $x_2y_1$ is not a signed simplicial edge. 
There is a negative edge $x_sy_t$ in the graph induced by $N(x_2y_1)$.
If $s = 1$ and $t = 2$, then $x_1,x_2,y_1,y_2$ induce the graph in $F_1$.
If $s = 1$ and $t \neq 2$ then $x_1,x_2,y_1,y_2,y_t$ induce a graph in $F_2$.
If $s \neq 1$ and $t = 2$, then $x_1,x_2,x_s,y_1,y_2$ also induce a graph in 
$F_2$. Finally, if $s \neq 1$ and $t \neq 2$ then $x_1,x_2,x_s,y_1,y_2,y_t$
induce a graph in $F_4$. These contradict the assumption. 
So $\widehat{G}$ contains no negative $P_4$.

Since no edge of $\widehat{G}$ is signed simplicial, there is a negative edge in 
the graph induced by $N(x_iy_j)$ for each $1 \leq i \leq \alpha$ and 
$1 \leq j \leq \beta$.
This implies that there are two negative edges which share no endvertex. 
Without loss of generality assume that $x_1y_1$ and $x_2y_2$ are 
negative. Let $x_sy_t$ be a negative edge in the graph induced by $N(x_1y_2)$. 
If $s \geq 3$ and $t \geq 3$, then $x_1,x_2,x_s,y_1,y_2,y_t$ induce 
a graph in $F_4$. Thus $s = 2$ or $t = 1$. Since $\widehat{G}$ contains no negative
$P_4$, $s \neq 2$ or $t \neq 1$. Hence either $s = 2$ and $t \neq 1$ or
$s \neq 2$ and $t = 1$. By symmetry we assume that $s = 2$ and $t \neq 1$. 
Let $x_py_q$ be a negative edge in the graph induced by $N(x_2y_1)$.
Since $\widehat{G}$ contains no negative $P_4$, $q \neq 2$ and $q \neq t$.
If $p \neq 1$ then $x_1, x_2, x_p, y_1, y_2, y_q$ induces a graph in $F_4$.
Thus $p = 1$ and $x_1,x_2,y_1,y_2,y_t,y_q$ induce a graph in $F_3$.
These again contradict the assumption. 
\qed

\begin{center}
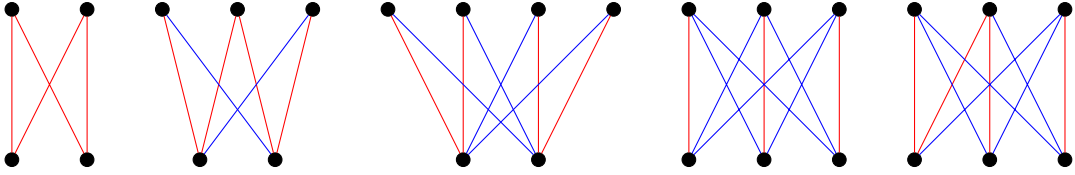
\begin{figure}[htb]
                \center
                \begin{tikzpicture}[>=latex]
                \begin{pgfonlayer}{nodelayer}
                %Fig a
\node [style=blackvertex] (1) at (0,0) {};
\node [style=blackvertex] (2) at (1,0) {};
\node [style=blackvertex] (3) at (0,2) {};
\node [style=blackvertex] (4) at (1,2) {};

\draw[style=edge] (1) to (3) to (2) to (4) to (1) [red];
%\node [style=textbox] at (.5,-.8) {$F^*_1$};

\node [style=blackvertex] (5) at (2.5,0) {};
\node [style=blackvertex] (6) at (3.5,0) {};
\node [style=blackvertex] (7) at (2,2) {};
\node [style=blackvertex] (8) at (3,2) {};
\node [style=blackvertex] (9) at (4,2) {};

\draw[style=edge] (7) to (5) to (8) to (6) to (9)[red];
\draw[style=edge] (6) to (7)[blue];
\draw[style=edge] (5) to (9)[blue];
%\node [style=textbox] at (3,-.8) {$F^*_2$};

\node [style=blackvertex] (10) at (6,0) {};
\node [style=blackvertex] (11) at (7,0) {};
\node [style=blackvertex] (12) at (5,2) {};
\node [style=blackvertex] (13) at (6,2) {};
\node [style=blackvertex] (14) at (7,2) {};
\node [style=blackvertex] (15) at (8,2) {};

\draw[style=edge] (12) to (10) to (13)[red];
\draw[style=edge] (14) to (11) to (15)[red];
\draw[style=edge] (12) to (11) to (13)[blue];
\draw[style=edge] (14) to (10) to (15)[blue];
%\node [style=textbox] at (6.5,-.8) {$F^*_3$};

\node [style=blackvertex] (16) at (9,0) {};
\node [style=blackvertex] (17) at (10,0) {};
\node [style=blackvertex] (18) at (11,0) {};
\node [style=blackvertex] (19) at (9,2) {};
\node [style=blackvertex] (20) at (10,2) {};
\node [style=blackvertex] (21) at (11,2) {};

\draw[style=edge] (16) to (19)[red];
\draw[style=edge] (17) to (20)[red];
\draw[style=edge] (18) to (21)[red];
\draw[style=edge] (16) to (20) to (18) to (19) to (17) to (21) to (16)[blue];
%\node [style=textbox] at (10,-.8) {$F^*_4$};

\node [style=blackvertex] (22) at (12,0) {};
\node [style=blackvertex] (23) at (13,0) {};
\node [style=blackvertex] (24) at (14,0) {};
\node [style=blackvertex] (25) at (12,2) {};
\node [style=blackvertex] (26) at (13,2) {};
\node [style=blackvertex] (27) at (14,2) {};

\draw[style=edge] (25) to (22) to (26) to (23)[red];
\draw[style=edge] (24) to (27)[red];
\draw[style=edge] (22) to (27) to (23) to (25) to (24) to (26)[blue];
%\node [style=textbox] at (13,-.8) {$F^{**}_4$};

                \end{pgfonlayer}
                \end{tikzpicture}
\caption{\label{minForb-complete}The minimal forbidden subgraphs for chordal
signed complete bigraphs}
        \end{figure}
\end{center}

Not all graphs in $F_1 \cup F_2 \cup F_3 \cup F_4$ are minimal non-chordal bigraphs. 
For instance, if any of the black edges in $F_2$ is replaced by a red edge then it 
contains the graph in $F_1$. On the other hand, replacing both black edges in $F_2$ 
by blue edges results in a minimal non-chordal bigraph. 
The five graphs in Figure \ref{minForb-complete} are all obtained in this way from 
the graphs in Figure \ref{Forb-complete}. By Theorem \ref{comp} they are not 
chordal. Deleting a vertex in any of them results in a chordal signed bigraph. 
So they are minimal non-chordal signed bigraphs. Moreover, every signed graph 
in $F_1 \cup F_2 \cup F_3 \cup F_4$ contains a graph in Figure \ref{minForb-complete}
as an induced subgraph. 
Hence by Theorem \ref{comp} if a signed complete bigraph is not chordal then 
it contains a graph in Figure \ref{minForb-complete} as an induced subgraph. 
Therefore we have the following:

\begin{theorem} \label{minF}
A signed complete graph is chordal if and only if it does not contain any 
graph in Figure \ref{minForb-complete} as an induced subgraph.
\qed
\end{theorem}

\section{Non-separable bigraphs} \label{3}

A bigraph is {\em separable} if it contains an induced $2K_2$, otherwise it is
{\em non-separable}. Since a cycle of length $\geq 6$ contains an induced $2K_2$, 
it cannot be an induced subgraph of a non-separable bigraph. Thus every non-separable
bigraph is a chordal bigraph. However, not all signed non-separable bigraphs are 
chordal signed bigraphs (e.g., the negative $C_4$ is non-separable but not a chordal 
signed bigraph). In this section we extend Theorem \ref{comp} to a forbidden 
subgraph characterization of chordal signed non-separable bigraphs.

We first take a look at some basic properties of non-separable bigraphs.
Let $G$ be a non-separable bigraph with bipartition $(X,Y)$. Since $G$ does not 
contain an induced $2K_2$, the neighbourhoods of the vertices in $X$ are
comparable and the neighbourhoods of the vertices of $Y$ are comparable.
It follows that the vertices of $X$ and $Y$ can be ordered 
$x_1, x_2, \dots, x_\alpha$ and $y_1, y_2, \dots, y_\beta$ respectively 
in such a way that
\[N(x_1) \supseteq N(x_2) \supseteq \cdots \supseteq N(x_\alpha)\ \mbox{and}\
N(y_1) \subseteq N(y_2) \subseteq \cdots \subseteq N(y_\beta).\] 
We call such a vertex ordering $x_1, x_2, \dots, x_\alpha, y_1, y_2, \dots, y_\beta$
a {\em canonical ordering} of $G$.

\begin{lemma} \label{basic_non-sep}
Let $G$ be a non-separable bigraph with bipartition $(X,Y)$ which has no isolated
vertex and let $x_1, x_2, \dots, x_\alpha, y_1, y_2, \dots, y_\beta$ be a 
canonical ordering of $G$. Then the following statements hold:
\begin{description}
\item{(1)} $x_1$ is adjacent to all vertices in $Y$.
\item{(2)} Each vertex in $N(y_1)$ is adjacent to all vertices in $Y$.
\item{(3)} For each $x_i \in N(y_1)$, $x_iy_1$ is a simplicial edge; in particular,
           $x_1y_1$ is a simplicial edge.
\item{(4)} If $x_iy_j$ is a simplicial edge then $G-x_iy_j$ is non-separable.
\item{(5)} If $x_iy_j$ is a non-simplicial edge then there exist
                   $x_k \in N(y_j)$ and $y_\ell \in N(x_i)$ with $k > i$ and 
                   $\ell < j$ such that $x_k$ and $y_\ell$ are not adjacent.  
\end{description}
\end{lemma}
\pf For each $y_j \in Y$, since $G$ has no isolated vertex, $y_j$ is adjacent to 
some $x_i \in X$. The canonical ordering ensures that $N(x_i) \subseteq N(x_1)$. 
Hence $y_j \in N(x_i) \subseteq N(x_1)$ and so $x_1$ is adjacent to $y_j$. 
This proves $(1)$.

For each $y_j \in Y$, the canonical ordering ensures that $N(y_1) \subseteq N(y_j)$.
Thus each vertex in $N(y_1)$ is adjacent to $y_j$. This proves $(2)$.

For each $x_i \in N(y_1)$, since $x_i$ is adjacent to all vertices in $Y$, $x_iy_1$ 
is an edge and $x_i$ is adjacent to all vertices in $N(x_1)$. So $x_iy_1$ is 
a simplicial edge and in particular $x_1y_1$ is a simplicial edge. This proves $(3)$.

To prove $(4)$, suppose that $G-x_iy_j$ contains an induced $2K_2$ consisting of 
edges $x_ay_b, x_cy_d$. Since $G$ is non-separable, the edges $x_ay_b, x_cy_d$ do not
form an induced $2K_2$ in $G$. So we must have $x_i = x_a$ and $y_j = y_d$ or 
$x_a = x_c$ and $y_j = y_b$. %as $G$ is non-separable. 
But in either case $x_iy_j$ is 
not a simplicial edge. 

To prove $(5)$, suppose that $x_iy_j$ is a non-simplicial edge. Then there 
are non-adjacent vertices $x_k \in N(y_j)$ and $y_\ell \in N(x_i)$. This implies that 
$N(x_i) \not\subseteq N(x_k)$ and $N(y_j) \not\subseteq N(y_\ell)$. The canonical
ordering ensures that $k > i$ and $\ell < j$.  
\qed

\begin{center}
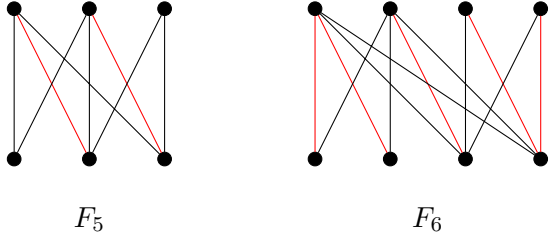
\begin{figure}[htb]
                \center
                \begin{tikzpicture}[>=latex]
                \begin{pgfonlayer}{nodelayer}
                %Fig a
\node [style=blackvertex] (1) at (0,0) {};
\node [style=blackvertex] (2) at (1,0) {};
\node [style=blackvertex] (3) at (2,0) {};
\node [style=blackvertex] (4) at (0,2) {};
\node [style=blackvertex] (5) at (1,2) {};
\node [style=blackvertex] (6) at (2,2) {};

\draw[style=edge] (2) to (4)[red];
\draw[style=edge] (3) to (5)[red];
\draw[style=edge] (1) to (4) to (3) to (6) to (2) to (5) to (1);
\node [style=textbox] at (1,-.8) {$F_5$};

\node [style=blackvertex] (7) at (4,0) {};
\node [style=blackvertex] (8) at (5,0) {};
\node [style=blackvertex] (9) at (6,0) {};
\node [style=blackvertex] (10) at (7,0) {};
\node [style=blackvertex] (11) at (4,2) {};
\node [style=blackvertex] (12) at (5,2) {};
\node [style=blackvertex] (13) at (6,2) {};
\node [style=blackvertex] (14) at (7,2) {};

\draw[style=edge] (7) to (11) to (8)[red];
\draw[style=edge] (13) to (10) to (14)[red];
\draw[style=edge] (9) to (12)[red];
\draw[style=edge] (7) to (12) to (10) to (11) to (9) to (13);
\draw[style=edge] (8) to (12);
\draw[style=edge] (9) to (14);
\node [style=textbox] at (5.5,-.8) {$F_6$};

                \end{pgfonlayer}
                \end{tikzpicture}
\caption{\label{Forb-non-sep}Additional forbidden subgraphs for chordal signed 
non-separable bigraphs}
        \end{figure}
\end{center}

Since complete bigraphs are non-separable, the graphs in 
$F_1 \cup F_2 \cup F_3 \cup F_4$ in Figure \ref{Forb-complete} are forbidden 
in chordal signed non-separable bigraphs. 
Two additional sets $F_5, F_6$ of forbidden subgraphs for chordal signed 
non-separable bigraphs are depicted in Figure \ref{Forb-non-sep}. 
We show that the graphs in $F_1 \cup F_2 \cup \cdots \cup F_6$ are the forbidden 
subgraphs for chordal signed non-separable bigraphs.

\begin{theorem}\label{non-separableForbidden}
Let $\widehat{G}$ be a signed non-separable bigraph. Then $\widehat{G}$ is chordal 
if and only if it does not contain any graph in $F_1 \cup F_2 \cup \cdots \cup F_6$ 
as an induced subgraph.
\end{theorem}

Before proving Theorem \ref{non-separableForbidden} we need two lemmas. 

\begin{lemma} \label{nonseparableZ1}
Let $\widehat{G}$ be a signed non-separable bigraph with bipartition $(X,Y)$ and
without isolated vertices and let
$x_1, x_2, \dots, x_\alpha, y_1, y_2, \dots, y_\beta$ be a canonical ordering of
$\widehat{G}$. Suppose that $\widehat{G}$ contains a graph in $Z_1$ induced by
$x_i, x_j, x_k, y_\ell, y_r$ (see in Figure \ref{Z1}) where $x_i, x_j \in N(y_1)$.
If no edge in $\widehat{G}$ is signed simplicial, then $\widehat{G}$ contains a graph
in $F_2 \cup F_3 \cup F_4 \cup F_5$ as an induced subgraph.
\end{lemma}

\begin{center}
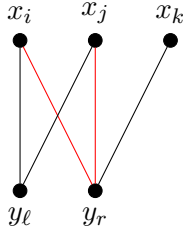
\begin{figure}[htb]
                \center
                \begin{tikzpicture}[>=latex]
                \begin{pgfonlayer}{nodelayer}
                %Fig a
\node [label={below:$y_\ell$}] [style=blackvertex] (1) at (0,0) {};
\node [label={below:$y_r$}] [style=blackvertex] (2) at (1,0) {};
\node [label={above:$x_i$}] [style=blackvertex] (3) at (0,2) {};
\node [label={above:$x_j$}] [style=blackvertex] (4) at (1,2) {};
\node [label={above:$x_k$}] [style=blackvertex] (5) at (2,2) {};

\draw[style=edge] (3) to (2) to (4)[red];
\draw[style=edge] (3) to (1) to (4);
\draw[style=edge] (2) to (5);
%\node [style=textbox] at (1,-.8) {$Z_1$};

                \end{pgfonlayer}
                \end{tikzpicture}
\caption{\label{Z1} Graph $Z_1$}
        \end{figure}
\end{center}

\pf Since $x_i, x_j \in N(y_1)$, $N(x_i)=N(x_j)=Y$ by Lemma \ref{basic_non-sep}(2).
Since $x_k$ is adjacent to $y_r$ but not to $y_\ell$, 
$N(y_r) \not\subseteq N(y_\ell)$.  
The canonical ordering ensures that $\ell < r$. 
Since $x_k$ is not adjacent to all vertices in $Y$,
$N(x_k) \subset N(x_i) = N(x_j) = Y$ and hence $k > i$ and $k > j$. 

Choose such a graph in $Z_1$ so that $k$ is as large as possible.  We claim that 
$x_ky_r$ is a simplicial edge in $G$. Indeed, if $x_ky_r$ is not a simplicial edge in
$G$. 
Then by Lemma \ref{basic_non-sep}(5) there exist $x_{k'} \in N(y_r)$ with $k' > k$ 
and $y_{\ell'} \in N(x_k)$ with $\ell' < r$ such that $x_{k'}$ and $y_{\ell'}$ are 
not adjacent. Thus
$x_i, x_j, x_{k'}, y_{\ell'}, y_r$ induce a graph in $Z_1$ with $k' > k$,
contradicting our choice. Hence $x_ky_r$ is a simplicial edge in $G$. 
By assumption $x_ky_r$ is not signed simplicial in $\widehat{G}$. So 
the subgraph of $\widehat{G}$ induced by $N(x_ky_r)$ has a negative edge. 

Let $x_ay_b$ be a negative edge in the subgraph induced by $N(x_ky_r)$. 
If $x_a$ is adjacent to $y_\ell$ then $\widehat{G}$ contains a graph in $F_5$ 
induced by
$x_i, x_a, x_k, y_\ell, y_r, y_b$ when $x_a \neq x_i$ or by
$x_j, x_a, x_k, y_\ell, y_r, y_b$ when $x_a \neq x_j$.
So we may assume that $x_a$ is not adjacent to $y_\ell$ which means that 
$x_a \neq x_i$ and $x_a \neq x_j$. Our choice of a graph in $Z_1$ implies that 
$a < k$.

Since $x_ay_r$ is not signed simplicial, the subgraph induced by $N(x_ay_r)$ is
either not a biclique or is a biclique that contains a negative edge.
Suppose first that the subgraph induced by $N(x_ay_r)$ is not a biclique. Then
there exist $x_c \in N(y_r)$ with $c > a$ and $y_d \in N(x_a)$ with $d <r$ such that 
$x_c$ and $y_d$ are not adjacent. We must have $c \leq k$ as otherwise 
$x_i, y_j, x_c, y_\ell, y_r$ would induce a graph in $Z_1$ with $c > k$, 
contradicting our choice. Thus $y_b \in N(x_k) \subseteq N(x_c)$. Hence 
$x_j, x_a, x_c, y_d, y_r, y_b$ induce a graph in $F_5$.
Suppose now that the subgraph induced by $N(x_ay_r)$ is a biclique that contains 
a negative edge $x_cy_d$. Note that $x_c$ is adjacent to $y_b$ as otherwise $c > k$ 
and $x_i, x_j, x_c, y_b, y_r$ induce a graph in $Z_1$, a contradiction to our choice.
If $y_d = y_b$ then $x_i, x_j, x_c, x_a, y_r, y_b$ induce a graph in $F_3$ 
when $x_c \neq x_i$ and $x_c \neq x_j$ and $x_i, x_j, x_a, y_r, y_b$ induce
a graph in $F_2$ when $x_c = x_i$ or $x_c = x_j$.
If $y_d \neq y_b$ %and $x_c$ is adjacent to $y_b$ 
then $\widehat{G}$ contains a graph in $F_4$ induced by $x_i, x_c, x_a, y_d, y_r, y_b$
when $x_c \neq x_i$ or by $x_j, x_c, x_a, y_d, y_r, y_b$ when $x_c \neq x_j$.
%If $y_d \neq y_b$ and $x_c$ is not adjacent to $y_b$ then $\widehat{G}$ contains 
%a graph in
 %$F_5$ induced by $x_i, x_c, x_a, y_d, y_r, y_b$ 
%when $x_c \neq x_i$ or by $x_j, x_c, x_a, y_d, y_r, y_b$ when $x_c \neq x_j$. 
\qed

\begin{center}
\begin{figure}[htb]
                \center
                \begin{tikzpicture}[>=latex]
                \begin{pgfonlayer}{nodelayer}
                %Fig a
\node [label={below:$y_\ell$}] [style=blackvertex] (1) at (0,0) {};
\node [label={below:$y_q$}] [style=blackvertex] (2) at (1,0) {};
\node [label={below:$y_r$}] [style=blackvertex] (3) at (2,0) {};
\node [label={above:$x_i$}] [style=blackvertex] (5) at (0,2) {};
\node [label={above:$x_j$}] [style=blackvertex] (4) at (1,2) {};
\node [label={above:$x_k$}] [style=blackvertex] (6) at (2,2) {};

\draw[style=edge] (1) to (5) to (2)[red];
\draw[style=edge] (4) to (3)[red];
\draw[style=edge] (1) to (4) to (2);
\draw[style=edge] (5) to (3) to (6);
%\node [style=textbox] at (1,-.8) {$Z_2$};

                \end{pgfonlayer}
                \end{tikzpicture}
\caption{\label{Z2} Graph $Z_2$}
        \end{figure}
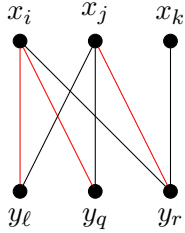
\end{center}

\begin{lemma} \label{nonseparableZ2}
Let $\widehat{G}$ be a signed non-separable bigraph with bipartition $(X,Y)$ and
without isolated vertices and let
$x_1, x_2, \dots, x_\alpha, y_1, y_2, \dots, y_\beta$ be a canonical ordering of
$\widehat{G}$. Suppose that $\widehat{G}$ contains a graph in $Z_2$ induced by
$x_i, x_j, x_k, y_\ell, y_q, y_r$ (see Figure \ref{Z2}) where $x_i, x_j \in N(y_1)$. 
If no edge in $\widehat{G}$ is signed simplicial, then $\widehat{G}$ contains 
a graph in $F_3 \cup F_4 \cup F_5 \cup F_6$ as an induced subgraph.
\end{lemma}
\pf Choose such a graph in $Z_2$ so that $k$ is as large as possible. 
We claim that the subgraph of $G$ induced by $N(x_ky_r)$ is a biclique. Indeed, if not
by Lemma \ref{basic_non-sep}(5), there exist non-adjacent $x_a \in N(y_r)$ with
$a > k$ and $y_b \in N(x_k)$ with $b < r$. Since $x_k$ is adjacent to $y_b$ but not
to $y_\ell$ or $y_q$, $N(y_b) \supset N(y_\ell)$ and $N(y_b) \supset N(y_q)$. 
This implies that $x_a$ is not adjacent to $y_\ell$ or $y_q$ as it is not adjacent
to $y_b$. Hence $x_i, x_j, x_a, y_\ell, y_q, y_r$ induce a graph in $Z_2$ with 
$a > k$, a contradiction to the choice of a graph in $Z_2$. So the subgraph induced 
by $N(x_ky_r)$ is a biclique. 
Since $x_ky_r$ is not signed simplicial in $\widehat{G}$, the biclique induced by 
$N(x_ky_r)$ contains a negative edge $x_ay_b$.

If $x_a = x_i$ then $x_i, x_j, x_k, y_\ell, y_b, y_r$ induce a graph in $F_5$.
If $x_a = y_j$ then $x_i, x_j, y_\ell, y_q, y_b, y_r$ induce a graph in $F_3$.
So assume $x_a \neq x_i$ and $x_a \neq x_j$. 
If $x_a$ is adjacent to $y_q$ then $x_j, x_a, x_k, y_q, y_b, y_r$ induce a graph in 
$F_5$.
If $x_a$ is adjacent to $y_l$ then $x_j, x_a, x_k, y_l, y_b, y_r$ induce a graph in
$F_5$.
Hence we further assume that $x_a$ is not adjacent to $y_q$ or $y_l$. 

Since $x_ay_r$ is not signed simplicial, the subgraph of $\widehat{G}$ induced by 
$N(x_ay_r)$ is either not a biclique or is a biclique that contains
a negative edge. Consider first the case when the subgraph induced by $N(x_ay_r)$ is 
not a biclique. By Lemma \ref{basic_non-sep}(5) there exist non-adjacent 
$x_c \in N(y_r)$ with $c > a$ and $y_d \in N(x_a)$ with $d < r$. 
The vertices $x_c$ and $y_b$ must be adjacent as otherwise 
$x_i, x_j, x_c, y_l, y_q, y_r$ induce a graph in $Z_2$ with $c > k$, contradicting
our choice. Then $x_j, x_a, x_c, y_d, y_b, y_r$ induce a graph in $F_5$.
%We see that $\widehat{G}$ contains a graph in $F_5$ induced by $x_j, x_a, x_c, y_d, y_b, y_r$. 
Consider now the case when the subgraph of $\widehat{G}$ induced by $N(x_ay_r)$ is 
a biclique which contains a negative edge $x_cy_d$. 
If $x_c = x_i$ then $\widehat{G}$ contains a graph in $F_5$ induced by 
$x_i, x_j, x_a, y_\ell, y_d, y_r$. 
If $x_c = x_j$ then $\widehat{G}$ contains a graph in $F_3$ induced by 
$x_i, x_j, y_\ell, y_q, y_d, y_r$. 
So assume that $x_c \neq x_i$ and $x_c \neq x_j$. 
If $x_cy_q$ is an edge then $\widehat{G}$ contains a graph in $F_5$ induced by 
$x_j, x_a, x_c, y_q, y_d, y_r$. 
If $x_cy_q$ is not an edge and $y_b \neq y_d$, then $\widehat{G}$ contains 
a graph in $F_4$ induced by $x_j x_a, x_c, y_b, y_d, y_r$. 
If $x_cy_q$ is not an edge and $y_b = y_d$ then $\widehat{G}$ 
contains a graph in $F_6$ induced by $x_i, x_j, x_a, x_c, y_\ell, y_q, y_b, y_r$.
\qed

\bigskip

{\bf Proof of Theorem \ref{non-separableForbidden}.} 
We only prove the if-part of the statement as the only-if part is discussed above. 
So assume that $\widehat{G}$ does not contain any graph in 
$F_1 \cup F_2 \cup \cdots \cup F_6$ as an induced subgraph. 
Suppose to the contrary that $\widehat{G}$ is not chordal. 
By Lemma \ref{basic_non-sep}(4) deleting a signed simplicial edge from $\widehat{G}$ 
maintains the property of being non-separable. Note also that deleting any signed
simplicial edge from $\widehat{G}$ does not results in a graph which contains a 
graph in $F_1 \cup F_2 \cup \cdots \cup F_6$. Hence we can further assume that 
$\widehat{G}$ has no signed 
simplicial edge and has no isolated vertex.

Let $(X,Y)$ be the bipartition and $x_1,x_2,\dots,x_\alpha,y_1,y_2,\dots,y_\beta$ be
a canonical ordering of $\widehat{G}$. Since no edge in $\widehat{G}$ is signed 
simplicial, for any edge $x_iy_j$ the subgraph induced by $N(x_iy_j)$ is not 
a positive biclique, that is, either it is not a biclique or is a biclique
but contains a negative edge. By Lemma \ref{basic_non-sep}(3), the subgraph induced 
by $N(x_1y_1)$ is a biclique so it contains a negative edge $x_ay_b$. 
Choose such an edge $x_ax_b$ so that $b$ is as large as possible. 
The vertices in $N(y_1)$ have the same neighbourhood by Lemma \ref{basic_non-sep}(2).
Thus by re-ordering the vertices of $N(y_1)$ if necessary we can assume that 
$x_a$ is the vertex with the largest subscript which is incident with
a negative edge in the subgraph induced by $N(x_1y_1)$. Hence the subgraph induced
by $N(x_1y_1)$ contains no negative edge $x_{a'}y_{b'}$ 
with $a' > a$ or $b' > b$. 

Consider the subgraph induced by $N(x_ay_1)$. It is a biclique according to
Lemma \ref{basic_non-sep}(3) and thus contains a negative edge $x_cy_d$. 
Note that $N(x_1) = Y$ by Lemma \ref{basic_non-sep}(1). So $x_cy_d$ is a negative
edge in the subgraph induced by $N(x_1y_1)$. The choice of $x_ay_b$ implies
that $1 \leq c < a$ and $1 < d \leq b$. 

{\bf Case~1.} $d < b$. 

If some vertex $x$ is adjacent to $y_d$ (and hence to $y_b$ as 
$N(y_d) \subseteq N(y_b)$) but not to a vertex $y \in Y$, then 
$x_c, x_a, x, y, y_d, y_b$ induce a graph in $F_5$, contradicting the assumption. 
So any vertex adjacent to $y_d$ is adjacent to all vertices in $Y$. This implies 
that the subgraph induced by $N(x_ay_d)$ is a biclique and hence contains 
a negative edge $x_ey_f$.
Since $x_e$ is adjacent to $y_d$, it is adjacent all vertices in $Y$ and in 
particular to $y_1$. So the edge $x_ey_f$ is a negative edge in the subgraph induced
by $N(x_1y_1)$. The choice of $x_ay_b$ implies $e < a$ and $f \leq b$.
If $x_e \neq x_c$ and $y_f \neq y_b$ then $x_e, x_c, x_a, y_f, y_d, y_b$ induce a 
graph in $F_4$, contradicting the assumption. Hence $x_e = x_c$ or $y_f = y_b$.

{\bf Subcase~1.1.} $x_e = x_c$.

Suppose that some vertex $x$ is adjacent to $y_b$ but not to $y_d$ (and hence not to 
$y_1$ as $N(y_1) \subseteq N(y_d)$).
If $y_f = y_b$ then $x_e, x, x_a, y_b, y_d$ induce a graph in $Z_1$. 
By Lemma \ref{nonseparableZ1}, $\widehat{G}$ contains a graph in 
$F_2 \cup F_3 \cup F_4 \cup F_5$ as an induced subgraph, a contradiction to 
assumption.
If $y_f \neq y_b$ then $x_e, x, x_a, y_b, y_d, y_f$ induce a graph in $Z_2$ when 
$x$ is not adjacent to $y_f$ and a graph in $F_5$ otherwise. This means that 
in view of Lemma \ref{nonseparableZ2} $\widehat{G}$ contains a graph in 
$F_3 \cup F_4 \cup F_5 \cup F_6$ as an induced subgraph, a contradiction.
So every vertex adjacent to $y_b$ is adjacent to $y_d$ and hence
to all vertices in $Y$. This implies that the subgraph induced by $N(x_cy_b)$ is 
a biclique and thus contains a negative edge $x_gy_h$. Since $x_g$ is adjacent to
$y_b$, it is adjacent to all vertices in $Y$. Thus $x_gy_h$ is in the subgraph 
induced by $N(x_1y_1)$. If $x_g \neq x_a$ and $y_h \neq y_d$, then 
$x_c, x_a, x_g, y_d, y_b, y_h$ induce a graph in $F_4$, a contradiction. 
If $x_g = x_a$ and $y_h \neq y_d$, then $x_c, x_a, y_d, y_b, y_h$ induce a graph in 
$F_2$
when $y_f = y_b$ and $x_c, x_a, x_g, y_f, y_h$ induce a graph in $F_3$ when $y_f \neq y_b$, 
contradicting the assumption.
If $x_g = x_a$ and $y_h = y_d$, then $x_c, x_a, y_d, y_b$ induce the graph in $F_1$ 
when $y_f = y_b$ and $x_c, x_a, y_d, y_f, y_b$ induce a graph in $F_2$ when 
$y_f \neq y_b$, 
contradicting the assumption. If $x_g \neq x_a$ and $y_h = y_d$, then 
$x_c, x_a, x_g, y_d, y_b$ induce a graph in $F_2$ when $y_f = y_b$ and
$x_c, x_a, x_g, y_d, y_f, y_b$ induce a graph in $F_4$ when $y_f \neq y_b$, again 
contradicting the assumption.

{\bf Subcase~1.2.} $x_e \neq x_c$.

We must have $y_f = y_b$. 
%as otherwise $x_c, x_e, x_a, y_d, y_f, y_b$ induce an $F_4$, a contradiction to the assumption. 
As in Subcase~1.1 we have that every vertex 
adjacent to $y_b$ is adjacent to $y_d$ and hence to all vertices in $Y$.
This together with Lemma \ref{basic_non-sep}(3) imply that the subgraph induced by 
$N(x_cy_b)$ is a biclique and thus contains a negative edge $x_iy_j$. 
Since $x_i$ is adjacent to $x_b$ it is adjacent to all vertices in $Y$. 
Thus $x_iy_j$ is in the subgraph induced by $N(x_1y_1)$. If $y_j \neq y_d$ 
then $x_c, x_e, x_i, y_d, y_j, y_b$ induce a graph in $F_4$ when $x_i \neq x_e$ and 
$x_c, x_i, x_a, y_d, y_j, y_b$ induce a graph in $F_4$ when $x_i \neq x_a$, contradicting 
the assumption. If $y_j = y_d$, then $x_c, x_e, x_i, x_a, y_d, y_b$ induce a graph in $F_3$ 
when $x_i \neq x_e$ or $x_a$ and $x_c, x_e, x_a, y_d, y_b$ induce a graph in $F_2$ when 
$x_i = x_e$ or $x_a$, contradicting the assumption.
 
{\bf Case~2.} $d = b$. 

If some vertex $x$ is adjacent to $y_b$ but not to a vertex $y \in Y$ then 
$x_c, x_a, x, y, y_b$ induce a graph in $Z_1$ and by Lemma \ref{nonseparableZ1} 
$\widehat{G}$ contains a graph in $F_2 \cup F_3 \cup F_4 \cup F_5$ as an induced 
subgraph, a contradiction to the assumption. 
So any vertex adjacent to $y_b$ is adjacent to all vertices of $Y$. 
This together with Lemma \ref{basic_non-sep}(3) imply that the subgraph induced by 
$N(x_ay_b)$ is a biclique and hence contains a negative edge $x_ey_f$. 
Note that $x_ey_f$ is an edge in the subgraph induced by $N(x_1y_1)$. 
The choice of $x_ay_b$ implies that $e < a$ and $f < b$.

{\bf Subcase 2.1.} $x_e = x_c$.

Consider the subgraph induced by $N(x_cy_b)$. By Lemma \ref{basic_non-sep}(3) and 
the fact that every vertex adjacent to $y_b$ is adjacent to all vertices in $Y$,    
it is a biclique and hence contains a negative edge $x_gy_h$. Again $x_gy_h$ is 
in the subgraph induced by $N(x_1y_1)$. If $x_g = x_a$ and $y_h = y_f$
then $x_c, x_a, y_f, y_b$ induce the graph in $F_1$. If $x_g \neq x_a$ and $y_h \neq y_f$ 
then $x_c, x_g, x_a, y_f, y_h, y_b$ induce a graph in $F_4$. If $x_g = x_a$ and 
$y_h \neq y_f$ then $x_c, x_a, y_f, y_h, y_b$ induce a graph in $F_2$. If $x_g \neq x_a$ and
$x_h = y_f$ then $x_c, x_g, x_a, y_f, y_d$ induce a graph in $F_2$. All these contradict
the assumption.

{\bf Subcase 2.2.} $x_e \neq x_c$.

Consider the subgraph induced by $N(x_ey_b)$. A similar argument as in Subcase 2.1
shows that it is a biclique and hence contains a negative edge $x_gy_h$ which is in 
the subgraph induced by $N(x_1y_1)$. If $y_h \neq y_f$ then $\widehat{G}$ 
contains a graph in $F_4$ induced by $x_e, x_g, x_a, y_f, y_h, y_b$ when $x_g \neq x_a$, and
by $x_c, x_g, x_e, y_f, y_b, y_b$ when $x_g \neq x_c$. If $y_h = y_f$ and 
$x_g \neq x_c$ or $x_a$ then $x_c, x_e, x_g, x_a, y_f, y_b$ induce a graph in $F_3$.
If $y_h = y_f$ and $x_g = x_c$ or $x_a$ then $\widehat{G}$ contains a graph in $F_2$ induced
by $x_c, x_e, x_a, y_f, y_b$.
\qed

\section{Separable bigraphs} \label{4}

A set $S$ of vertices in a bigraph $G$ is {\em separating} if $G-S$ contains at 
least two non-trivial components. Clearly, a bigraph has a separating set
if and only if it is separable. A separating set is {\em minimal} if no proper 
subset of the set is separating. Suppose that $S$ is a minimal separating set. 
The minimality of $S$ ensures that each vertex in $S$ has a neighbour in every 
non-trivial component of $G-S$. In particular, there are non-trivial components $H$ 
and $H'$ in $G-S$ such that each vertex in $S$ has neighbours in $H$ and in $H'$.

In general, for a set $S$ of vertices in a bigraph $G$ and two non-trivial components
$H$ and $H'$ in $G-S$, we say that $S$ {\em minimally separates} $H$ and $H'$ if 
each vertex in $S$ has neighbours in $H$ and in $H'$. From above we know that
if $G$ is separable then there is a set $S$ of vertices in $G$ which minimally 
separates two non-trivial components in $G-S$.  

Our goal in this section is to find all forbidden subgraphs for chordal signed 
separable bigraphs. Signed cycles in $C_{2k}$ with $k \geq 3$ (see 
Figure~\ref{Forb-sep1}) are separable but not chordal, and deleting any vertex 
results in a chordal signed bigraph. So they are forbidden subgraphs for 
chordal signed separable bigraphs. The signed bigraphs in $D$ as depicted in 
Figure~\ref{Forb-sep1} are another set of such bigraphs.

\begin{center}
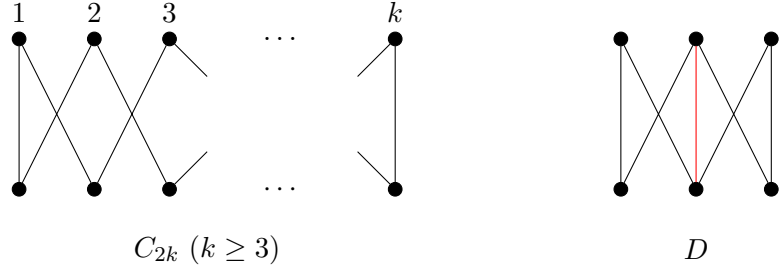
\begin{figure}[htb]
                \center
                \begin{tikzpicture}[>=latex]
                \begin{pgfonlayer}{nodelayer}
                %Fig a
\node [style=blackvertex] (1) at (0,0) {};
\node [style=blackvertex] (2) at (1,0) {};
\node [style=blackvertex] (3) at (2,0) {};
\node [style=blackvertex] (4) at (5,0) {};
\node [label={above:$1$}] [style=blackvertex] (5) at (0,2) {};
\node [label={above:$2$}] [style=blackvertex] (6) at (1,2) {};
\node [label={above:$3$}] [style=blackvertex] (7) at (2,2) {};
\node [label={above:$k$}] [style=blackvertex] (8) at (5,2) {};

\draw[style=edge] (3) to (6) to (1) to (5) to (2) to (7);
\draw[style=edge] (4) to (8);
\draw[style=edge] (3) to (2.5,.5);
\draw[style=edge] (7) to (2.5,1.5);
\draw[style=edge] (4) to (4.5,.5);
\draw[style=edge] (8) to (4.5,1.5);
\node [style=textbox] at (2.5,-.8) {$C_{2k}$ ($k \geq 3$)};
\node [style=textbox] at (3.5,0) {$\cdots$};
\node [style=textbox] at (3.5,2) {$\dots$};

\node [style=blackvertex] (9) at (8,0) {};
\node [style=blackvertex] (10) at (9,0) {};
\node [style=blackvertex] (11) at (10,0) {};
\node [style=blackvertex] (12) at (8,2) {};
\node [style=blackvertex] (13) at (9,2) {};
\node [style=blackvertex] (14) at (10,2) {};

\draw[style=edge] (10) to (13)[red];
\draw[style=edge] (9) to (13) to (11) to (14) to (10) to (12) to (9);
\node [style=textbox] at (9,-.8) {$D$};

                \end{pgfonlayer}
                \end{tikzpicture}
\caption{\label{Forb-sep1}Forbidden subgraphs $C_{2k}$ and $D$ for chordal signed 
separable bigraphs}
        \end{figure}
\end{center}

Denote ${\cal C} = \bigcup_{k\geq 3} C_{2k}$, that is, ${\cal C}$ consists of signed 
cycles of length $\geq 6$.

\begin{lemma} \label{min-ss} 
Let $\widehat{G}$ be a signed separable bigraph which does not contain a signed graph
in ${\cal C} \cup D$ in Figure \ref{Forb-sep1} as an induced subgraph. 
Suppose that $S$ minimally separates $H$ and $H'$. Then 
\begin{itemize}
\item $S$ induces a positive biclique in $\widehat{G}$, and 
\item any two vertices in $S$ from the same partite set of $\widehat{G}$ have 
a common neighbour in $H$ and a common neighbour in $H'$.
\end{itemize}
\end{lemma} 
\pf Let $u, u'$ be two vertices in $S$.
Since $S$ minimallly separates $H$ and $H'$, each of $u, u'$ has a neighbour in $H$ 
and a neighbour in $H'$. Let $P$ (respectively, $Q$) be a shortest path in $H$ 
(respectively, $H'$) from a vertex in $N(u)$ to a vertex in $N(u')$. 
If $u, u'$ are from different partite sets of $\widehat{G}$, then $uu'$ is an edge 
and $P, Q$ each has length one, as otherwise at least one of $uPu'Qu$, $uPu'u$, and 
$u'Quu'$ is an induced cycle of length $\geq 6$ in $\widehat{G}$, contradicting 
the assumption. Since $\widehat{G}$ does not contain a graph in $D$ as an induced 
subgraph, the edge $uu'$ is positive. Hence $S$ induces a positive biclique in 
$\widehat{G}$. If $u, u'$ are from the same partite set of $\widehat{G}$, then 
$uPu'Qu$ is an induced cycle of length $\geq 6$, unless $P$ and $Q$ each has length 0,
that is, $u, u'$ have a common neighbour in $H$ and a common neighbour in $H'$. 
\qed

\begin{center}
\begin{figure}[htb]
                \center
                \begin{tikzpicture}[>=latex]
                \begin{pgfonlayer}{nodelayer}
\node [style=blackvertex] (1) at (-2,0) {};
\node [label={right:$x$}] [style=blackvertex] (2) at (0,0) {};
\node [style=blackvertex] (3) at (-1,1) {};
\node [style=blackvertex] (4) at (-1,-1) {};

\draw[style=edge] (3) to (2) to (4)[red];
\draw[style=edge] (3) to (1) to (4);
\node [style=textbox] at (-1,-1.5) {$W_1$};

\node [style=blackvertex] (1) at (3,0) {};
\node [style=blackvertex] (2) at (5,0) {};
\node [label={right:$x$}] [style=blackvertex] (3) at (6,0) {};
\node [style=blackvertex] (4) at (4,1) {};
\node [style=blackvertex] (5) at (4,-1) {};

\draw[style=edge] (5) to (2) to (4)[red];
\draw[style=edge] (5) to (1) to (4);
\draw[style=edge] (2) to (3);
\draw[style=edge] (1) to [out=90,in=180] (4.5,1.5) to [out=0,in=90] (3)[red]; 
\node [style=textbox] at (4,-1.5) {$W_3$};

\node [style=blackvertex] (1) at (9,0) {};
\node [style=blackvertex] (2) at (11,0) {};
\node [label={right:$x$}] [style=blackvertex] (3) at (12,0) {};
\node [label={right:$x'$}] [style=blackvertex] (4) at (10,1) {};
\node [style=blackvertex] (5) at (10,-1) {};

\draw[style=edge] (2) to (4)[red];
\draw[style=edge] (4) to (1) to (5) to (2) to (3);
\draw[style=edge] (1) to [out=90,in=180] (10.5,1.5) to [out=0,in=90] (3)[red];
\node [style=textbox] at (10,-1.5) {$W_5$};

\node [style=blackvertex] (1) at (-2,-4) {};
\node [style=blackvertex] (2) at (0,-4) {};
\node [label={right:$x$}] [style=blackvertex] (3) at (1,-4) {};
\node [style=blackvertex] (4) at (-1,-3) {};
\node [style=blackvertex] (5) at (-1,-5) {};

\draw[style=edge] (5) to (2) to (4)[red];
\draw[style=edge] (5) to (1) to (4);
\draw[style=edge] (2) to (3);
%\draw[style=edge] (1) to [out=90,in=180] (4.5,1.5) to [out=0,in=90] (3)[red];
\node [style=textbox] at (-1,-5.5) {$W_2$};

\node [style=blackvertex] (1) at (3,-4) {};
\node [style=blackvertex] (2) at (5,-4) {};
\node [style=blackvertex] (3) at (6,-4) {};
\node [style=blackvertex] (4) at (4,-3) {};
\node [style=blackvertex] (5) at (4,-5) {};
\node [label={right:$x$}] [style=blackvertex] (6) at (7,-4){};

\draw[style=edge] (5) to (2) to (4)[red];
\draw[style=edge] (4) to (1) to (5);
\draw[style=edge] (2) to (3) to (6);
\draw[style=edge] (1) to [out=90,in=180] (4.5,-2.5) to [out=0,in=90] (3)[red];
\node [style=textbox] at (4,-5.5) {$W_4$};

\node [style=blackvertex] (1) at (9,-4) {};
\node [style=blackvertex] (2) at (11,-4) {};
\node [label={right:$x$}] [style=blackvertex] (3) at (12,-4) {};
\node [label={right:$x'$}] [style=blackvertex] (4) at (10,-3) {};
\node [style=blackvertex] (5) at (10,-5) {};

\draw[style=edge] (2) to (4)[red];
\draw[style=edge] (3) to (2) to (5) to (1) to (4);
%\draw[style=edge] (1) to [out=90,in=180] (10.5,-2.5) to [out=0,in=90] (3)[red];
\node [style=textbox] at (10,-5.5) {$W_6$};

                \end{pgfonlayer}
                \end{tikzpicture}
\caption{\label{Ws} Graphs $W_i$ ($1 \leq i \leq 6$)}
        \end{figure}
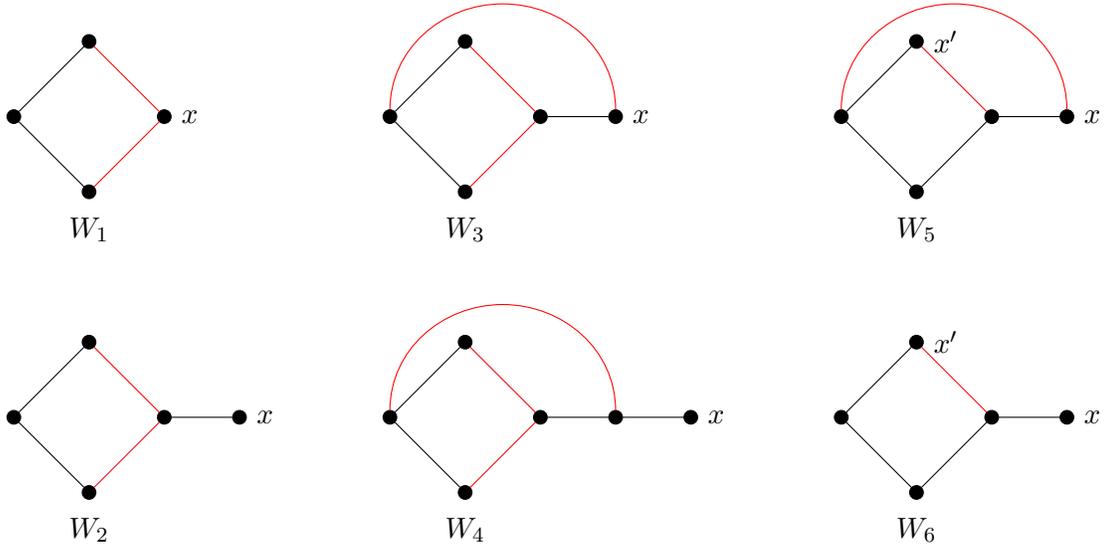
\end{center}

\begin{lemma} \label{s1}
Let $\widehat{G}$ be a chordal signed non-separable bigraph with bipartition $(X,Y)$.
Suppose that $S$ is a set of vertices such that
\begin{itemize}
\item $S$ induces a positive biclique in $\widehat{G}$,
\item every vertex of $S$ has a neighbour in $\widehat{G} - S$, and
\item $\widehat{G} - S$ is connected and contains at least one edge but none of them 
      is a signed simplicial edge of $\widehat{G}$.
\end{itemize}
Then $\widehat{G}$ contains a graph in $W_1 \cup W_2 \cup \cdots \cup W_6$ 
in Figure~\ref{Ws} as an 
induced subgraph where vertices $x,x'$ are the only ones in $S$ (and rest are in 
$\widehat{G} - S$).  
In particular, when $\widehat{G}$ is a signed complete bigraph, it contains a graph 
in $W_1 \cup W_3 \cup W_5$ as an induced subgraph.
\end{lemma}
\pf Since $\widehat{G}$ is a chordal signed bigraph, it does not contain any graph 
in $F_1 \cup F_2 \cup \cdots \cup F_6$ (see Figures \ref{Forb-complete} and 
\ref{Forb-non-sep}) as an induced subgraph. This fact will be used throughout 
the proof.

Consider first the case when $G$ is complete. We claim that there is a negative 
edge between $S$ and $\widehat{G} - S$. Indeed, since $\widehat{G}$ contains edges 
which are not chordal, there is at least one negative edge. There is no negative edge
with its both endvertices in $S$ as $S$ induces a positive biclique in $\widehat{G}$. 
By assumption $\widehat{G}$ is chordal and $\widehat{G} - S$ contains no signed 
simplicial edge of $\widehat{G}$. Thus any signed simplicial edge of $\widehat{G}$
has at least one endvertex in $S$. Let $uv$ be a signed simplicial edge of 
$\widehat{G}$ with $u \in S$. As $G$ is complete $\widehat{G} - \{u,v\}$ (which
is the subgraph of $\widehat{G}$ induced by $N(uv)$) is a positive biclique. 
Thus any negative edge is incident with $u$ or $v$. If there is no negative edge
between $S$ and $\widehat{G}-S$ then $v$ is in $\widehat{G}-S$ and all negative edges
of $\widehat{G}$ are in $\widehat{G}-S$ and incident with $v$.  We see now that 
any negative edge is in $\widehat{G} - S$ and is signed simplicial in $\widehat{G}$, 
a contradiction to the assumption. Hence there is a negative edge between $S$ and 
$\widehat{G} - S$.

Let $x \in S$ be a vertex incident with a negative edge. Suppose that $x$ is incident
with two negative edges, that is, there are vertices $y', y''$ in $\widehat{G} - S$ 
such that $xy', xy''$ are both negative. Since $\widehat{G} - S$ is connected it 
contains a vertex $x'$ is the opposite partite set of $y', y''$. 
Then $x, x', y', y''$ induce a graph in $W_1$ in $\widehat{G}$
with $x$ being the only vertex in $S$.
So assume that no vertex in $S$ incident with two negative edges. 
If some vertex $y \in S$ in the opposite partite set of $x$ which is adjacent to 
a negative edge, then the two negative edges incident with $x$ and $y$ must be 
adjacent to any signed simplicial edge of $\widehat{G}$. Thus $xy$ is 
signed simplicial. Together with the fact each 
of $x, y$ is incident with exactly one negative edge we conclude that the 
negative edges incident with $x, y$ are the only negative edges in $\widehat{G}$. 
This means that the edge in $\widehat{G} - S$ adjacent to the two negative edges 
is a signed simplicial edge in $\widehat{G}$, contradicting the assumption.
Hence no vertex of $S$ in the opposite partite set of $x$ is incident with a 
negative edge. Suppose that $x$ is the only vertex in $S$ that is incident with 
a negative edge. Denote this edge by $xy$. Let $x'$ be any vertex in $\widehat{G}-S$ 
and in the same partite set as $x$. Since $x'y$ is not a signed simplicial edge 
there is a negative edge $x''y''$ in the graph induced by $N(x'y)$. Then $xy''$ is a
signed simplicial edge and thus all negative edges other than $xy$ are incident with
$y''$. Since $x''y$ lies in $\widehat{G} - S$, it is not a signed simplicial edge in
$\widehat{G}$. There must be a vertex $x'''$ (which may be $x'$ but not $x''$) such
that $x'''y''$ is negative. We see that $x, x'', x''', y, y''$ induce a graph in
$W_3$ in $\widehat{G}$ with $x$ being the only vertex in $S$.
Finally suppose that $x$ is not the only vertex in $S$ that is incident with a 
negative edge. If all negative edges between $S$ and $\widehat{G} - S$ are incident 
with the same vertex $y$ then $y$ must be incident with any signed simplicial edge. 
It follows that $y$ is incident with all negative edges. But then any edge of 
$\widehat{G} - S$ incident with $y$ is a signed simplicial edge in $\widehat{G}$,
a contradiction. Hence there vertices $x, x' \in S$ incident with negative edges
$xy, x'y'$ where $y \neq y'$. Therefore $x, y, x', y'$ together with any vertex of 
$\widehat{G} - S$ in the same partite set as $x$ induce a graph in $W_5$ 
with $x, x'$ being the only vertices in $S$.  

Consider next the case when $G$ is non-separable but not complete. Let 
$$\prec:\ x_1, x_2, \dots, x_\alpha, y_1, y_2, \dots, y_\beta$$
be a canonical ordering of $G$. Since $G$ is not complete, $x_{\alpha}$ and $y_1$ 
are not adjacent and hence they are not both in $S$. By renaming the vertices if 
necessary we assume that $y_1$ is in $\widehat{G}-S$. Since $\widehat{G}-S$ 
is connected and has at least one edge, it contains edge $x_ky_1$ for some $k$. 
Choose such an edge with $k$ as small as possible. 
The canonical ordering $\prec$ ensures that $N(x_i) \supseteq N(x_k)$ for all 
$i < k$. Thus $x_iy_1$ is an edge for all $i < k$. 
The choice of $k$ implies that $x_i \in S$ for all $i < k$. The edge $x_ky_1$ is 
a simplicial edge in $G$ according to Lemma \ref{basic_non-sep}(3) but 
is not a signed simplicial edge in $\widehat{G}$ by assumption. 
So the subgraph of $\widehat{G}$ induced by $N(x_ky_1)$ contains a negative edge 
$x_ay_b$.

{\bf Case~1.} $x_a \in S$ and $y_b \notin S$.

The edge $x_ky_b$ is an edge in $\widehat{G} - S$ and by assumption is not signed
simplicial. That is, the subgraph of $\widehat{G}$ induced by $N(x_ky_b)$ is 
either not a biclique or is a biclique and contains a negative edge. 

{\bf Subcase 1.1.} The subgraph induced by $N(x_ky_b)$ is a biclique. 

Let $x_cy_d$ be a negative edge where $x_c \in N(y_b)$ and $y_d \in N(x_k)$. 
Since the subgraph induced by $N(x_ky_b)$ is a biclique and $x_a \in N(y_b)$, 
$N(x_k) \subseteq N(x_a)$ which implies that $x_ay_d$ is an edge.
Suppose that $x_c \in S$. Since $S$ induces a positive biclique and $x_cy_d$ is 
a negative edge, $y_d \notin S$. When $x_c = x_a$, $\widehat{G}$ contains a graph in 
$W_1$ 
induced by $x_a, x_k, y_b, y_d$ with $x_a$ being the only vertex in $S$; when 
$x_c \neq x_a$, $\widehat{G}$ contains a graph in $W_5$ induced by 
$x_a, x_c, x_k, y_b, y_d$ with $x_a, x_c$ being the only vertices in $S$. 

Suppose now that $x_c \notin S$. Since the subgraph induced by $N(x_ky_b)$ is a 
biclique, $x_c \in N(y_b)$ and $y_1 \in N(x_k)$, we have $N(x_k) \subseteq N(x_c)$ 
and $x_cy_1$ is an edge. The choice of $x_k$ thus implies that $k < c$. 
Moreover, the canonical ordering $\prec$ implies that $N(x_c) \subseteq N(x_k)$.
The edge $x_cy_b$ is in $\widehat{G}-S$ and by assumption is not signed simplicial. 
So the subgraph graph induced by $N(x_cy_b)$ is either not a biclique or a biclique 
and contains a negative edge. That is, there exist $x_e \in N(y_b)$ and 
$y_f \in N(x_c)$ such that $x_ey_f$ is either not an edge or is a negative edge. 
We claim that $x_ey_f$ is an edge (and hence a negative edge). 
Indeed, using the fact $N(x_ky_b)$ is a biclique, we have 
$y_f \in N(x_c) \subseteq N(x_k) \subseteq N(x_e)$, which $x_ey_f$ is an edge 
and hence a negative edge. Note that $x_ay_f$ is also an edge as 
$y_f \in N(x_k) \subseteq N(x_a)$. 

If $x_e = x_a$ then $y_f \notin S$ and $\widehat{G}$ contains a graph in $W_1$ induced by 
$x_a, x_k, y_b, y_f$ with $x_a$ being the only vertex in $S$. 
So assume that $x_e \neq x_a$. Suppose that $x_e = x_k$. If $y_f \neq y_d$ then
$x_a, x_c, x_e, y_b, y_d, y_f$ induce a graph in $F_4$, a contradiction.
So $y_f = y_d$. If $y_d \in S$ then $\widehat{G}$ contains a graph in $W_1$ induced by
$x_c, x_e, y_d, y_b$ with $y_d$ being the only vertex in $S$; 
if $y_d \notin S$ then $\widehat{G}$ contains a graph in $W_3$ induced by
$x_a, x_c, x_e, y_b, y_d$ with $x_a$ being the only vertex in $S$. 

Suppose that $x_e \neq x_k$. If $y_f \neq y_d$, then $x_a, x_c, x_e, y_b, y_d, y_f$
induce a graph in $F_4$, a contradiction. So $y_f = y_d$. 
If $y_d \in S$ then $\widehat{G}$ contains a graph in
$W_1$ induced by $x_c, x_e, y_1, y_d$ with $y_d$ being the only vertex in $S$; 
if $y_d \notin S$ and $x_e \notin S$ then $\widehat{G}$ contains a graph in $W_3$ induced by 
$x_a, x_c, x_e, y_b, y_d$ with $x_a$ being the only vertex in $S$; if $y_d \notin S$ 
and $x_e \in S$ then $\widehat{G}$ contains a graph in $W_5$ induced by $x_a, x_k, x_e, y_b, y_d$
with $x_a, x_e$ being the only vertices in $S$. 
   
{\bf Subcase 1.2.} The subgraph induced by $N(x_ky_b)$ is not a biclique.

Let $x_c \in N(y_b)$ and $y_d \in N(x_k)$ be chosen so that $x_c, y_d$ are not
adjacent and $c$ is as large as possible. 
By Lemma \ref{basic_non-sep}(5), $c > k$ and $d < b$. 
Since $y_d \in N(x_k) \subseteq N(x_a)$, $x_ay_d$ is an edge. The canonical ordering
$\prec$ ensures that $c > a$.% and $c > k$.

If $x_c \in S$ then $\widehat{G}$ contains a graph in $W_6$ induced by $x_a, x_k, x_c, y_1, y_b$ 
with $x_a, x_c$ being the only vertices in $S$. So assume that $x_c \notin S$. 
Then the edge $x_cy_b$ is in $\widehat{G} - S$ and by assumption is not signed 
simplicial. If there exist $x_e \in N(y_b)$ and $y_f \in N(x_c)$ such that 
$x_e, y_f$ are not adjacent then $e > c$ and $f < b$ by Lemma \ref{basic_non-sep}(5). 
Since $y_f \in N(x_c) \subseteq N(x_k)$, $x_ky_f$ is an edge. Thus $x_e, y_f$ is 
a pair of non-adjacent vertices with $x_e \in N(y_b)$, $y_f \in N(x_k)$ and $e > c$,
which contradicts to choice of the pair $x_c, y_d$. Hence the subgraph induced by
$N(x_cy_b)$ is a biclique and contains a negative edge. 

Let $x_ey_f$ be a negative edge where $x_e \in N(y_b)$ and $y_f \in N(x_c)$.  
Suppose that $x_ey_d$ is an edge. If $x_e = x_a$ then $\widehat{G}$ contains a graph in $W_1$ 
induced by $x_a, x_k, y_b, y_f$ with $x_a$ being the only vertex in $S$; 
if $x_e \neq x_a$ then $x_a, x_c, x_e, y_b, y_d, y_f$ induce a graph in $F_5$,
a contradiction. So ssume now that $x_ey_d$ is not an edge. 
We must have $e < c$ by the choice of the pair $x_c, y_d$. If $x_e \in S$ then 
$\widehat{G}$ contains a graph in $W_6$ induced by $x_a, x_k, x_e, y_1, y_b$ with 
$x_a, x_e$ 
being the only vertices in $S$. So assume that $x_e \notin S$. Then the edge $x_ey_b$ 
is in $\widehat{G}-S$ and by assumption is not signed simplicial. 
Then there exist $x_g \in N(y_b)$ and $y_h \in N(x_e)$ such that $x_gy_h$ is either
not an edge or a negative edge.

Suppose that $x_gy_h$ is not an edge. Note that $x_ky_h$ is an edge as 
$y_h \in N(x_e) \subseteq N(x_k)$. By the choice of $x_c, y_d$, $g \leq c$. Note also
that $x_gy_f$ is an edge (as $N(x_cy_b)$ is a biclique). So $y_f \neq y_h$. Then 
$x_a, x_g, x_e, y_b, y_f, y_h$ induce a graph in $F_5$, a contradiction.

Suppose that $x_gy_h$ is a negative edge. Note that $x_gy_f$ is an edge 
(as $N(x_cy_b)$ is a biclique). If $x_g = x_a$ then $\widehat{G}$ contains a graph in $W_1$ 
induced by $x_a, x_k, y_b, y_h$ with $x_a$ being the only vertex in $S$. So assume 
that $x_g \neq x_a$. If $y_h \neq y_f$ then $x_a, x_e, x_g, y_b, y_f, y_h$ induce
a graph in $F_4$, a contradiction. Thus $y_h = y_f$. If $y_f \in S$ then 
$\widehat{G}$ contains a graph in $W_1$ induced by $x_g, x_e, y_b, y_f$ with
$y_f$ being the only vertex in $S$; if $y_f \notin S$ and $x_g \notin S$ then 
$\widehat{G}$ contains a graph in $W_3$ induced by $x_a, x_e, x_g, y_b, y_f$ with 
$x_a$ being the only vertex in $S$; if $y_f \notin S$ and $x_g \in S$, then 
$\widehat{G}$ contains a graph in $W_5$ induced by $x_a, x_c, x_g, y_b, y_f$ with 
$x_a, x_g$ being the only vertices in $S$.  

{\bf Case~2.} $x_a \notin S$ and $y_b \in S$.

Then $x_ay_1$ is an edge in $\widehat{G}-S$ and thus not a signed simplicial edge. 
By Lemma \ref{basic_non-sep}(3), $N(x_ay_1)$ is a biclique in $G$. Hence there is 
negative edge $x_cy_d$ with $x_c \in N(y_1)$ and $y_d \in N(x_a)$. 
If $y_b = y_d$ then $\widehat{G}$ contains a graph in $W_1$ induced by $x_a, x_c, y_1, y_b$
with $y_b$ being the only vertex in $S$. So assume that $y_b \neq y_d$. 
If $y_d \in S$ then $x_c \notin S$ and $\widehat{G}$ contains a graph in $W_5$ induced by 
$x_a, x_c, y_1, y_b, y_d$ with $y_b, y_d$ being the only vertices in $S$. Hence 
we can assume from now on that $y_d \notin S$. This means that $x_ay_d$ is in 
$\widehat{G} - S$ and by assumption is not a signed simplicial edge. 
Then there exist $x_e \in N(y_d)$ and $y_f \in N(x_a)$ such that either $x_ey_f$ 
is not an edge or a negative edge when the subgraph induced by $N(x_ay_d)$ is 
a biclique.

Suppose that $x_ey_f$ is a negative edge in it. If $y_f = y_b$ then $x_e \notin S$ as 
$x_ey_b$ is a negative edge. We see that $\widehat{G}$ contains a graph in $W_1$ induced by 
$x_a, x_e, y_1, y_b$ with $y_b$ being the only vertex in $S$. So assume that 
$y_f \neq y_b$. If $x_c \neq x_e$ then $x_a, x_c, x_e, y_b, y_d, y_f$ induce a graph 
in $F_4$, a contradiction. Hence $x_c = x_e$. If $x_c \in S$ then $\widehat{G}$ 
contains a graph in $W_1$ induced by $x_a, x_c, y_d, y_f$ with $x_c$ being the only 
vertex in $S$; if $x_c \notin S$ and $y_f \notin S$ then $\widehat{G}$ contains 
a graph in $W_3$ induced by $x_a, x_c, y_b, y_d, y_f$ with $y_b$ being
the only vertex in $S$; if $x_c \notin S$ and $y_f \in S$ (note $y_f \neq y_1$) then
$\widehat{G}$ contains a graph in $W_5$ induced by $x_a, x_c, y_1, y_b, y_f$ with 
$y_b, y_f$ being the only vertices in $S$.  

Suppose now $x_ey_f$ is not an edge. We may assume such a pair of non-adjacent 
vertices $x_e, y_f$ in $N(x_ay_d)$ is chosen so that $e$ is the largest. 
Note that $x_ey_b$ cannot be an edge as otherwise $x_a, x_c, x_e, y_b, y_d, y_f$ 
induce a graph in $F_5$, is a contradiction. 
Since $x_e, y_b$ are not adjacent and $y_b \in S$, $x_e \notin S$.
Then $x_ey_d$ is in $\widehat{G}-S$ and not a signed simplicial edge. 
We claim that $N(x_ey_d)$ is a biclique. Indeed, for any $x_g \in N(y_d)$, 
we must have $g < e$ as otherwise $x_g, y_f$ would be a pair in $N(x_ay_d)$ with
$g > e$, contradicting the choice of $x_e, y_f$. Combining this with 
Lemma \ref{basic_non-sep}(5) we conclude that $N(x_ey_d)$ is a biclique containing
a negative edge $x_gy_h$.

Suppose that $x_gy_b$ is an edge. If $x_g \neq x_c$ then 
$x_c, x_e, x_g, y_b, y_d, y_h$ induce a graph in $F_5$, a contradiction. 
So $x_g = x_c$.
If $x_c \in S$ then $\widehat{G}$ contains a graph in $W_1$ induced by 
$x_a, x_c, y_d, y_h$. Assume that $x_c \notin S$.
If $y_h \notin S$ then $\widehat{G}$ contains a graph in $W_2$ induced by 
$x_c, x_e, y_b, y_d, y_h$ with $y_b$ being the only vertex in $S$;
if $y_h \in S$ then $\widehat{G}$ contains a graph in $W_5$ induced by 
$x_a, x_c, y_1, y_b, y_h$ with $y_b, y_h$ being the only vertices in $S$.

Suppose $x_g, y_b$ are not adjacent. Then $x_g \notin S$ (as $y_b \in S$). Thus
$x_gy_d$ is in $\widehat{G}-S$ and is not a signed simplicial edge.  
Let $x_i, y_j$ be a pair of vertices in $N(x_gy_d)$ such that either they are not
adjacent or forming a negative edge.
Suppose that $x_i = x_a$. 
Note that $y_j \neq y_b$ as $x_gy_b$ is not an edge but $x_gy_j$ is an edge.
Suppose that $x_i = x_a$.
Then $x_ay_j$ is a negative edge and $x_a, x_c, x_g, y_b, y_d, y_j$ induce
a graph in $F_5$, a contradiction.
%If $y_j = y_h$ then $\widehat{G}$ contains $F_5$ induced by 
%$x_a, x_c, x_g, y_b, y_d, y_j$, a contradiction; if $y_j \neq y_h$ then 
%$\widehat{G}$ contains $F_4$ induced by $x_a, x_c, x_g, y_d, y_h, y_j$, also
%a contradiction. 
So $x_i \neq x_a$. 

Suppose that $x_i = x_c$. If $x_c \in S$ then $\widehat{G}$ contains a graph in $W_1$ induced 
by $x_a, x_c, y_d, y_j$ with $x_c$ being the only vertex in $S$; 
if $x_c \notin S$ and $y_j \notin S$ then $\widehat{G}$ contains a graph in $W_2$ induced 
by $x_g, x_c, y_b, y_d, y_j$ with $y_b$ being the only vertex in $S$; 
if $x_c \notin S$ and $y_j \in S$ then $\widehat{G}$ contains  a graph in $W_5$ induced by
$x_a, x_c, y_b, y_1, y_j$ with $y_b, y_j$ being the only vertices in $S$.
Hence we can assume that $x_i \neq x_c$. 

Suppose that $x_iy_b$ is an edge. 
If $y_j \neq y_b$ then $x_a, x_c, x_i, y_b, y_d, y_j$ induce a graph in $F_4$ or 
$x_a, x_c, x_i, y_b, y_d, y_j$ induce a graph in $F_5$ (depending on whether or not 
$x_i, y_j$ are adjacent), a contradiction. So $y_j = y_b$ and 
$\widehat{G}$ contains a graph in $W_1$ induced
by $x_a, x_i, y_b, y_d$ with $y_b$ being the only vertex in $S$.  
So we can assume that $x_i, y_b$ are not adjacent. 

Suppose that $x_i \neq x_e$. The choice of $x_e, y_f$ implies that $i < e$ 
(as otherwise $x_i, y_f$ would have been chosen). Note that $x_iy_h$ is an edge as 
$y_h \in N(x_e) \subseteq N(x_i)$. Hence $x_c, x_g, x_i, y_d, y_h, y_j$ induce 
a graph in $F_4$ or a graph in $F_5$ (depending on whether or not $x_i, y_j$ are 
adjacent), a contradiction.

Suppose that $x_i = x_e$. If $y_j \neq y_h$ then $x_c, x_i, x_g, y_d, y_j, y_h$ 
induce a graph in $F_4$ or a graph in $F_5$ (depending on whether or not 
$x_i, y_j$ are adjacent), a contradiction. So $y_j = y_h$ and $\widehat{G}$ contains 
a graph in $W_4$ induced by $x_c, x_g, x_i, y_b, y_d, y_h$ with $y_b$ being the only 
vertex in $S$.  

{\bf Case~3.} $x_a \notin S$ and $y_b \notin S$.

By Lemma \ref{basic_non-sep}(3), the subgraph induced by $N(x_{k'}y_1)$ is a biclique
in $G$ for any $k'$ (including the case when $k' = k$). In view of Cases 1 and 2, we 
can assume that no negative edge in this biclique has an endvertex in $S$. 

Since $x_ay_1$ is in $\widehat{G}-S$, it is not signed simplicial. 
Since the subgraph induced by $N(x_ay_1)$ is a biclique in $G$, it contains 
a negative edge $x_cy_d$. By our assumption neither $x_c$ nor $y_d$ is in $S$. 

{\bf Subcase 3.1} $x_c = x_k$.

Since $x_ky_b$ is not signed simplicial, there exist $x_e \in N(y_b)$ and 
$y_f \in N(x_k)$ such that either $x_ey_f$ is not an edge or a negative edge when
the subgraph induced by $N(x_ky_b)$ is a biclique. 

Suppose that the subgraph induced by $N(x_ky_b)$ is a biclique and thus $x_ey_f$
is a negative edge. Since $x_e \in N(y_1)$, by Lemma \ref{basic_non-sep}(2), $x_e$ 
is adjacent to every vertex in $Y$. By our assumption above, $x_e$ is not in $S$.  
Assume first that $x_e \neq x_a$ and $y_b = y_d$.  

Since $x_ey_b$ is an edge in $\widehat{G} -S$ and hence not signed simplicial,
there exist $x_g \in N(y_b)$ and $y_h \in N(x_e)$ such that either $x_gy_h$ 
is not an edge or a negative edge in the biclique induced by $N(x_ey_b)$.
We claim that the subgraph induced by $N(x_ey_b)$ is not a biclique. Indeed,
if the subgraph induced by $N(x_ey_b)$ is a biclique, then we have the following:
\begin{itemize}
\item when $y_h = y_f$ and $x_g \in \{x_a, x_k\}$,
$x_a, x_c, x_e, y_f, y_b$ induce a graph in $F_2$;
\item when $y_h = y_f$ and $x_g \notin \{x_a, x_k\}$,
$x_a, x_c, x_e, x_g, y_b, y_f$ induce a graph $F_3$;
\item when $y_h \neq y_f$ and $x_g \in \{x_a, x_k\}$,
$x_a, x_c, x_e, y_b, y_f, y_h$ induce a graph in $F_4$;
\item when $y_h \neq y_f$ and $x_g \notin \{x_a, x_k\}$, 
$x_a, x_e, x_g, y_b, y_f, y_h$ induce a graph in $F_4$.
\end{itemize}
So $x_gy_h$ is not an edge.

Denote $A = N(y_d) \setminus N(y_h)$. Then $A \neq \emptyset$ as $x_g \in A$.
If $S$ contains a vertex $x_t \in A$ then $\widehat{G}$ contains $W_2$ induced by
$x_k, x_a, x_t, y_1, y_d$ with $x_t$ being the only vertex in $S$. Otherwise
$S$ contains no vertex of $A$. That is, for any $x_t \in A$, $x_ty_d$ is not a
signed simplicial edge. Then $\widehat{G}$ contains  a graph in $Z_1$ induced by
$x_k, x_a, x_t, y_1, y_d$.  Applying Lemma \ref{nonseparableZ1} to the subgraph of
$\widehat{G}$ induced by $A \cup \{x_k, x_a, y_1, y_d\}$, we conclude that
$\widehat{G}$ contains a graph in $F_2 \cup F_3 \cup F_4 \cup F_5$ as an induced
subgraph, a contradiction.

Assume now that $x_e \neq x_a$ and $y_b \neq y_d$. Then $y_f = y_d$
as otherwise $x_a, x_e, x_k, y_b, y_d, y_f$ induce a graph in $F_4$,
a contradiction. Since $x_ay_d$ is in $\widehat{G}-S$ it is not signed simplicial. 
We claim that the subgraph induced by $N(x_ay_d)$ is not a biclique. Indeed, if it is
then it contains a negative edge $x_gy_h$. If $y_h \neq y_b$ and $x_k = x_g$ then 
$x_a, x_e, x_g, y_b, y_f, y_h$ induce a graph in $F_4$;
if $y_h \neq y_b$ and $x_k \neq x_g$ then $x_a, x_k, x_g, y_b, y_d, y_h$ a graph in 
$F_4$; if $y_h = y_b$, and $x_g = x_k$ (or $x_g = x_e$), 
then $x_a, x_k, x_e, y_b, y_d$ induce a graph in $F_2$; 
if $y_h = y_b$ and $x_g \neq x_k$ (and $x_g \neq x_e$) then 
$x_a, x_k, x_e, x_g, y_b, y_d$ induce a graph in $F_3$, contradictions.
So there are non-adjacent vertices $x_g \in N(y_d)$ and $y_h \in  N(x_a)$.

If $S$ contains a vertex $x_t \in A$ then $\widehat{G}$ contains  a graph in $W_2$ induced by
$x_k, x_e, x_t, y_1, y_d$ with $x_t$ being the only vertex in $S$. Otherwise 
$S$ contains no vertex of $A$. That is, for any $x_t \in A$, $x_ty_d$ is not a
signed simplicial edge. Then $\widehat{G}$ contains a graph in $Z_1$ induced by 
$x_k, x_e, x_t, y_1, y_d$.  Applying Lemma \ref{nonseparableZ1} to the subgraph of 
$\widehat{G}$ induced by $A \cup \{x_k, x_e, y_1, y_d\}$, we conclude that
$\widehat{G}$ contains a graph in $F_2 \cup F_3 \cup F_4 \cup F_5$ as an induced
subgraph, a contradiction.

Assume now that $x_e = x_a$. Note that $x_ay_d$ is not signed simplicial.
We claim that the subgraph induced by $N(x_ay_d)$ is not a biclique. Indeed, if it
is then it contains a negative edge $x_gy_h$. If $y_b = y_d$, $x_g = x_k$ and
$y_f = y_h$, then $x_a, x_k, y_b, y_h$ induce the graph in $F_1$; 
if $y_b = y_d$, $x_g = x_k$ and $y_f \neq y_h$, then $x_a, x_k, y_b, y_f, y_h$
induce a graph in $F_2$; if $y_b = y_d$, $x_g \neq x_k$ and $y_h = y_f$, then 
$x_a, x_g, x_k, y_b, y_f$ induce a graph in $F_2$; 
if $y_b = y_d$, $x_g \neq x_k$ and $y_h \neq y_f$, then 
$x_a, x_k, x_g, y_b, y_f, y_h$ induce a graph in $F_4$;
if $y_b \neq y_d$ and $x_g \neq x_k$, then c$x_a, x_k, x_g, y_b, y_f, y_h$
induce a graph in $F_4$;
if $y_b \neq y_d$, $x_g = x_k$ and $y_h = y_b$ (or $y_h = y_f$), then 
$x_a, x_k, y_b, y_d, y_f$ induce a graph in $F_2$;
if $y_b \neq y_d$, $x_g = x_k$ and $y_h \neq  y_b$ (and $y_h \neq y_f$), then 
$x_a, x_k, y_b, y_d, y_f, y_h$ induce a graph in $F_3$.
So there exist non-adjacent vertices $x_g \in N(y_d)$ and $y_h \in N(x_a)$. 

Suppose that $y_b = y_d$. A similar proof as above (with the same definition of 
the set $A$) shows that $\widehat{G}$ contains a graph in $W_2$ induced by 
$x_a, x_c, x_t, y_1, y_d$ when $S$ contains a vertex $x_t \in A$ with $x_t$ being 
the only vertex in $S$; otherwise, $\widehat{G}$ contains a graph in  
$F_2 \cup F_3 \cup F_4 \cup F_5$ as an induced subgraph, 
a contradiction.

Suppose that $y_b \neq y_d$. Then $x_g$ is adjacent to neither of $y_b, y_f$ as 
otherwise $x_a, x_k, x_g, y_b, y_d, y_h$ or by $x_a, x_k, x_g, y_d, y_f, y_h$ would
induce a graph in $F_5$,
If some vertex $x_t \in N(y_d) \setminus N(y_h)$ is adjacent to
$y_f$ then $x_a, x_k, x_t, y_d, y_f, y_h$ induce a graph in $F_5$.
If some vertex $x_t \in N(y_d) \setminus N(y_h)$ is adjacent to $y_f$ or $y_b$ then 
$x_a, x_k, x_t, y_d, y_f, y_h$ or by $x_a, x_k, x_t, y_d, y_b, y_h$ would induce 
a graph in $F_5$.
So assume that none of vertices in $N(y_d) \setminus N(y_h)$ is adjacent to $y_f$ or 
$y_b$. 
%So assume that none of vertices in $N(y_d) \setminus N(y_h)$ is adjacent to $y_f$.
Again, a similar proof as above (with the same definition of the set $A$ but 
using Lemma \ref{nonseparableZ2}) shows that $\widehat{G}$ contains a graph in $W_4$ 
induced 
by $x_a, x_k, x_t, y_b, y_d, y_f$ when $S$ contains $x_t \in A$ with $x_t$ being 
the only vertex in $S$; otherwise, $\widehat{G}$ contains a graph in 
$F_3 \cup F_4 \cup F_5 \cup F_6$ as an induced subgraph,
a contradiction. Here we apply Lemma \ref{nonseparableZ2} (instead of 
Lemma \ref{nonseparableZ1}) to the subgraph of $\widehat{G}$ induced by 
$A \cup \{x_a, x_k, y_b, y_d, y_f\}$ and verify that it contains  a graph in $Z_2$.

Suppose now that $x_ey_f$ is not an edge. Since $x_ay_1$ is an edge the canonical 
ordering implies that $x_ay_f$ is an edge. It follows that $x_e \neq x_a$.

Suppose that $y_b = y_d$. We must have $x_e \in S$ as otherwise a similar proof as 
above with the set $A = N(y_d) \setminus N(y_f)$ shows that  $\widehat{G}$ contains a graph in
$Z_1$ induced by $x_a, x_c, x_e, y_f, y_b$. 
By Lemma \ref{nonseparableZ1} $\widehat{G}$ contains a graph in 
$F_2 \cup F_3 \cup F_4 \cup F_5$ as an induced subgraph, 
a contradiction. Hence $\widehat{G}$ contains a graph in $W_2$ induced by 
$x_a, x_k, x_e, y_1, y_b$ with $x_e$ being the only vertex in $S$.

Suppose now that $y_b \neq y_d$. Then $x_ey_d$ is not an edge as otherwise
$x_a, x_k, x_e, y_b, y_d, y_f$ induce a graph in $F_5$, a contradiction.
Since $x_ay_d$ is an edge in $\widehat{G} -S$, it is not signed simplicial,
Then there exist $x_g \in N(y_d)$ and $y_h \in N(x_a)$ such that either $x_gy_h$ is
not an edge or a negative edge when the subgraph induced by $N(x_ay_d)$ is a biclique.
Suppose that $x_gy_h$ is not an edge. Note that $x_gy_d$ is an edge but $x_ey_d$
is not so, by the canonical ordering, $g < e$ and hence that $x_gy_b$ is an edge but
$x_ey_h$ is not an edge. Clearly, $x_g \neq x_a$ or $x_k$.
We see that $x_a, x_k, x_g, y_b, y_d, y_h$ induce a graph in $F_5$,
a contradiction. Hence the subgraph induced by $N(x_ay_d)$ is a biclique
(and thus $x_gy_h$ is a negative edge). Since $x_g \in N(y_d)$ and $y_1 \in N(x_a)$,
$x_gy_1$ is an edge. The canonical ordering then implies $x_g$ is adjacent to every
vertex in $Y$. In the case when $x_g = x_k$, we must have that $x_e, y_h$ are not
adjacent as otherwise $x_a, x_k, x_e, y_b, y_d, y_h$ induce a graph in $F_5$,
a contradiction. Moreover, we must also have $x_e \in S$, as otherwise 
$\widehat{G}$ contains a graph in $Z_2$ induced by
$x_a, x_k, x_e, y_b, y_d, y_h$. (Here we use a similar proof as above with
$A = N(y_b) \setminus N(y_f)$.) By Lemma \ref{nonseparableZ2} 
$\widehat{G}$ contains a graph in $F_3 \cup F_4 \cup F_5 \cup F_6$ as 
an induced subgraph, again a contradiction.
We see now that $\widehat{G}$ contains a graph in $W_4$ induced by 
$x_a, x_e, x_k, y_b, y_d, y_h$ with $x_e$ being the only vertex in $S$. 
In the case when $x_g \neq x_k$, we must have $y_h = y_b$ as otherwise
$\widehat{G}$ contains a graph in $F_4$ induced by $x_a, x_k, x_g, y_b, y_d, y_h$,
a contradiction. Further, we must have $x_e \in S$ as otherwise $\widehat{G}$ 
contains a graph in $Z_1$ induced by $x_a, x_g, x_e, y_f, y_b$, and by Lemma \ref{nonseparableZ1}
$\widehat{G}$ contains a graph in $F_2 \cup F_3 \cup F_4 \cup F_5$ as 
an induced subgraph, again a contradiction.
Therefore $\widehat{G}$ contains a graph in $W_2$ induced by $x_a, x_e, x_g, y_1, y_b$ with
$x_e$ being the only vertex in $S$.

{\bf Subcase 3.2.} $x_c \neq x_k$.

We know from above that $x_ay_d$ is not signed simplicial. Thus there exist
$x_e \in N(y_d)$ and $y_f \in N(x_a)$ such that either $x_ey_f$ is not an edge or
is a negative edge in the case when the subgraph induced by $N(x_ay_d)$ is a biclique.

Suppose that the latter case occurs and thus $x_ey_f$ is a negative edge. 
So $x_e, y_f \notin S$ according to the assumption at the beginning of Case~3.
Since the subgraph induced by $N(x_ay_d)$ is a biclique and $x_e \in N(y_d)$, 
$N(x_a) \subseteq N(x_e)$ which implies that $x_e$ is adjacent to every vertex in 
$Y$. 

Note that $x_ey_b$ is not a signed simplicial edge as it is in $\widehat{G}-S$. 
There exist $x_g \in N(y_b)$ and $y_h \in N(x_e)$ such 
that either $x_gy_h$ is not an edge or a negative edge in the biclique induced 
by $N(x_ey_b)$. 

Assume that $y_b = y_d$.
We show by contradiction that the subgraph induced by $N(x_ey_b)$ is not a biclique.
So assume that the subgraph induced by $N(x_ey_b)$ is a biclique and $x_gy_h$ is 
a negative edge. Clearly, $x_g$ is adjacent to every vertex in $Y$.
If $y_h \neq  y_f$ and $x_e \neq x_c$ then $\widehat{G}$ contains a graph in $F_4$ 
induced by
$x_a, x_e, x_g, y_b, y_f, y_h$ when $x_g = x_c$ or by 
$x_c, x_e, x_g, y_b, y_f, y_h$ when $x_g \neq x_c$;
if $y_h \neq y_f$ and $x_e = x_c$ then $\widehat{G}$ contains a graph in $F_2$ 
induced by 
$x_a, x_c, y_f, y_b, y_h$ when $x_g = x_a$ or contains a graph in $F_4$ induced by 
$x_a, x_e, x_g, y_b, y_f, y_h$ when $x_g \neq x_a$;
if $y_h = y_f$ and $x_e = x_c$ then $\widehat{G}$ contains the graph in $F_1$ 
induced by 
$x_a, x_c, y_b, y_f$ when $x_g = x_a$ or contains a graph in $F_2$ induced by 
$x_a, x_c, x_g, y_b, y_f$ when $x_g \neq x_a$;
if $y_h = y_f$ and $x_e \neq x_c$ then $\widehat{G}$ contains a graph in $F_2$ 
induced by 
$x_a, x_c, x_e, y_b, y_f$ when $x_g = x_a$ or $x_g = x_c$, or 
contains a graph in $F_3$ induced by 
$x_a, x_c, x_e, x_g, y_b, y_f$ when $x_g \neq x_a$ and $x_g \neq x_c$.  
Hence $x_gy_h$ is not an edge. Clearly, $x_g \notin \{x_a, x_c, x_k\}$.
Therefore $\widehat{G}$ contains a graph in $W_2$ induced by $x_a, x_c, x_t, y_1, y_b$ with 
$x_t$ being the only vertex in $S$. 

Assume now that $y_b \neq y_d$. We separate this in two cases. First we consider
the case when $x_c = x_e$. Since $x_ey_b$ is not signed simplicial, 
there exist $x_g \in N(y_b)$ and $y_h \in N(x_e)$ such that either $x_gy_h$ is not
an edge or is a negative edge in the biclique induced by $N(x_ey_b)$. 
We claim that the subgraph induced by $N(x_ey_b)$ is not a biclique. Suppose not.
Then $x_gy_h$ is a negative edge. 
If $y_f = y_b$ and $y_h = y_d$ then $\widehat{G}$ contains the graph in $F_1$ 
induced by $x_a, x_c, y_b, y_d$ when $x_g = x_a$ or a graph in $F_2$ induced by 
$x_a. x_e, x_g, y_b, y_d$ when $x_g \neq x_a$;
if $y_f = y_b$ and $y_h \neq y_d$ then $\widehat{G}$ contains a graph in $F_2$ 
induced by $x_a, x_c, y_b, y_d, y_h$ when $x_g = x_a$ or in $F_4$ induced by
$x_a, x_c, x_g, y_b, y_d, y_h$ when $x_g \neq x_a$;
if $y_f \neq y_b$ and $y_h = y_d$ then $\widehat{G}$ contains a graph in $F_2$ 
induced by $x_a, x_c, y_b, y_f, y_h$ when $x_g = x_a$ or in $F_4$ induced by
$x_a, x_e, x_g, y_b, y_f, y_h$ when $x_g \neq x_a$;
if $y_f \neq y_b$ and $y_h \neq y_d$ then $\widehat{G}$ contains a graph in $F_4$ 
induced by
$x_a, x_c, x_g, y_b, y_d, y_h$ when $x_g \neq x_a$ or in $F_3$ induced by 
$x_a, x_c, y_b, y_d, y_f, y_h$ when $x_g = x_a$ and $y_f \neq y_h$ or in
$F_2$ induced by $x_a, x_c, y_b, y_d, y_f$ when $x_g = x_a$ and $y_f = y_h$.
Hence the subgraph induced by $N(x_ey_b)$ is not a biclique, that is, 
$x_g, y_h$ are not adjacent.

We claim that $x_g$ is in $S$ and adjacent to neither of $y_d, y_f$.
Indeed, $\widehat{G}$ contains a graph in $F_5$ induced by 
$x_a, x_c, x_g, y_b, y_d, y_h$ when $x_g$ is adjacent to $y_d$ or induced by
$x_a, x_c, x_g, y_b, y_f, y_h$ when $x_g$ is adjacent to $y_f$, contradicting 
the assumption. If $x_g \notin S$ then a similar proof as above (with 
$A = N(y_b) \setminus (N(y_d) \cup N(y_f))$) shows that $\widehat{G}$ contains 
a graph in $Z_2$ and hence by
Lemma \ref{nonseparableZ2} a graph in $F_3 \cup F_4 \cup F_5 \cup F_6$ as 
an induced subgraph, a contradiction.  
Therefore $\widehat{G}$ contains a graph in $W_4$ induced by 
$x_a, x_c, x_g, y_b, y_d, y_f$ with $x_g$ being the only vertex in $S$. 

We consider now the case when $x_c \neq x_e$. We must have $y_f = y_b$ as otherwise
$\widehat{G}$ contains a graph in $F_4$ induced by $x_a, x_c, x_e, y_b, y_d, y_f$, 
a contradiction. Note that $x_cy_b$ is not a signed simplicial edge as it is in
$\widehat{G} - S$. There exist $x_g \in N(y_b)$ and $y_h \in N(x_c)$ such that
either $x_gy_h$ is not an edge or a negative edge in the biclique induced by
$N(x_cy_b)$. 

We claim that the subgraph induced by $N(x_cy_b)$ is not a biclique. Suppose that it
is. Then $x_gy_h$ is a negative edge. If $y_h = y_d$ then $\widehat{G}$ contains 
a graph in $F_2$ induced by  
$x_a, x_c, x_e, y_b, y_d$ when $x_g = x_a$ or $x_g = x_e$ or in $F_3$ induced by 
$x_a, x_c, x_e, x_g, y_b, y_d$ when $x_g \neq x_a$ and $x_g \neq x_e$;
if $y_h \neq y_d$ then $\widehat{G}$ contains a graph in $F_4$ induced by $x_c, x_e, x_g, y_b, y_d, y_h$ when $x_g = x_a$ or induced by $x_a, x_c, x_g, y_b, y_d, y_h$ when 
$x_g \neq x_a$.
Hence the subgraph induced by $N(x_cy_b)$  is not a biclique and $x_g, y_h$ are 
not adjacent.
  
We must have $x_g \in S$, otherwise $\widehat{G}$ contains a graph in $Z_1$ and hence a graph in $F_2 \cup F_3 \cup F_4 \cup F_5$ as an induced subgraph, a contradiction. Therefore $\widehat{G}$ contains a graph in $W_2$ induced by
$x_a, x_e, x_g, y_1, y_b$ with $x_g$ being the only vertex in $S$.

Suppose now that $x_e, y_f$ are not adjacent. We separate this in two cases. 
First assume that $y_b = y_d$. Then $x_e \in S$ as otherwise $\widehat{G}$ contains a graph in 
$Z_1$ and hence a graph in $F_2 \cup F_3 \cup F_4 \cup F_5$ as an induced subgraph, a contradiction. 
(Here again we use a similar proof with $A = N(y_d) \setminus N(y_f)$.) 
Therefore $\widehat{G}$ contains a graph in 
$W_2$ induced by $x_a, x_c, x_e, y_1, y_b$ with $x_e$ being the only vertex in $S$.
Assume now that $y_b \neq y_d$. We must have that $x_e, y_b$ are not adjacent as 
otherwise $\widehat{G}$ contains a graph in $F_5$ induced by $x_a, x_c, x_e, y_b, y_d, y_f$,
a contradiction. Since $x_cy_b$ is not a signed simplicial edge, there exist
$x_g \in N(y_b)$ and $y_h \in N(x_c)$ such that either $x_gy_h$ is not an edge or 
a negative edge in the biclique induced by $N(x_cy_b)$.  
We claim that $N(x_cy_b)$ induces a biclique. Suppose not; $x_g, y_h$ are not 
adjacent. Note $g < e$ because $y_b \in N(x_g) \setminus N(x_e)$, which implies
$x_gy_d$ is an edge. But then $\widehat{G}$ contains a graph in $F_5$ induced by 
$x_a, x_c, x_g, y_b, y_d, y_h$, a contradiction.   
Therefore $x_gy_h$ is a negative edge. Note that $x_g$ is adjacent to every vertex in
$Y$.  

If $y_h = y_d$ then $x_e \in S$ as otherwise $\widehat{G}$ contains a graph in $Z_1$ 
and hence a graph in $F_2 \cup F_3 \cup F_4 \cup F_5$ as an induced subgraph, a contradiction. Hence $\widehat{G}$ contains a graph in $W_2$ induced by
$x_c, x_e, x_g, y_1, y_d$ with $x_e$ being the only vertex in $S$.
So assume $y_h \neq y_d$. If $x_ey_h$ is an edge then $\widehat{G}$ contains a graph in $F_5$
induced by 
$x_c, x_e, x_g, y_b, y_d, y_h$, a contradiction. So $x_e, y_h$ are not adjacent.
Then we must have $x_e \in S$ as otherwise $\widehat{G}$ contains a graph in $Z_2$ and hence 
a graph in $F_3 \cup F_4 \cup F_5 \cup F_6$ as an induced subgraph, a contradiction. Therefore $\widehat{G}$ contains a graph in $W_4$ induced by
$x_a, x_c, x_e, y_b, y_d, y_h$ with $x_e$ being the only vertex in $S$.
This completes the proof.
\qed 

\begin{lemma} \label{h1orh2}
Let $\widehat{G}$ be a chordal signed separable bigraph. Suppose that $S$ minimally
separates $\widehat{H_1}$ and $\widehat{H_2}$. Then 
$\widehat{H_1}$ or $\widehat{H_2}$ must contain a signed simplicial edge of
$\widehat{G}$. 
\end{lemma}
\pf Let $\widehat{G'}$ be the subgraph of $\widehat{G}$ induced by 
$S \cup V(\widehat{H_1}) \cup V(\widehat{H_2})$. Since $\widehat{G}$ is chordal, 
$\widehat{G'}$ is also chordal and hence contains a signed simplicial edge $e$. 
As each vertex of $S$ has a neighbour 
in $\widehat{H_1}$ and a neighbour in $\widehat{H_2}$, neither of the endvertices of
$e$ is in $S$. Since no edge has one endvertex in $\widehat{H_1}$ and the other 
in $\widehat{H_2}$, $e$ is an edge of $\widehat{H_1}$ or of $\widehat{H_2}$. Hence 
the neighbourhood of $e$ in $\widehat{G'}$ is equal to the neighbourhood of $e$ in
$\widehat{G}$, which means that $e$ is a signed simplicial edge of $\widehat{G}$.
\qed

\begin{center}
\begin{figure}[htb]
                \center
                \begin{tikzpicture}[>=latex]
                \begin{pgfonlayer}{nodelayer}
\node [style=blackvertex] (1) at (-1,0) {};
\node [style=blackvertex] (2) at (1,0) {};
\node [style=blackvertex] (3) at (2,0) {};
\node [style=blackvertex] (4) at (0,1) {};
\node [style=blackvertex] (5) at (0,-1) {};
\node [style=blackvertex] (6) at (3,0) {};
\node [style=textbox] at (4,0) {$\dots$};
\node [label={below:$x$}] [style=blackvertex] (7) at (5,0) {};

\node [style=textbox] at (2.5,-1.5) {${\cal T}_1$};

\draw[style=edge] (5) to (2) to (4)[red];
\draw[style=edge] (5) to (1) to (4);
\draw[style=edge] (2) to (3.5,0);
\draw[style=edge] (7) to (4.5,0);
%\draw[style=edge] (1) to [out=90,in=180] (4.5,1.5) to [out=0,in=90] (3)[red];

\node [style=blackvertex] (1) at (7,0) {};
\node [style=blackvertex] (2) at (9,0) {};
\node [style=blackvertex] (3) at (10,0) {};
\node [style=blackvertex] (4) at (8,1) {};
\node [style=blackvertex] (5) at (8,-1) {};
\node [style=blackvertex] (6) at (11,0) {};
\node [style=textbox] at (12,0) {$\dots$};
\node [label={below:$x$}] [style=blackvertex] (7) at (13,0) {};

\node [style=textbox] at (10.5,-1.5) {${\cal T}_2$};

\draw[style=edge] (5) to (2) to (4)[red];
\draw[style=edge] (5) to (1) to (4);
\draw[style=edge] (2) to (11.5,0);
\draw[style=edge] (7) to (12.5,0);
\draw[style=edge] (1) to [out=90,in=180] (8.5,1.5) to [out=0,in=90] (3)[red];
                \end{pgfonlayer}
                \end{tikzpicture}
\caption{\label{Ts} Tadpoles}
        \end{figure}
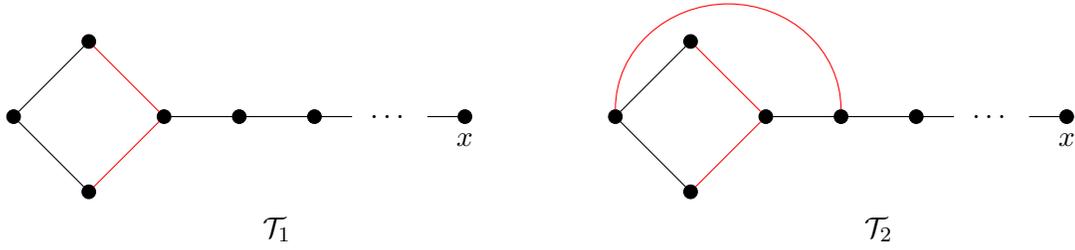
\end{center}

Figure \ref{Ts} depicts two sets ${\cal T}_1, {\cal T}_2$ of signed graphs.
We shall call the graphs in ${\cal T}_1 \cup {\cal T}_2$ {\em tadpoles}, and refer to each graph in ${\cal T}_1$ as a tadpole of {\em type 1} and each graph in 
${\cal T}_2$ a tadpole of {\em type 2}. We call the vertex $x$ the {\em end} of 
the tadpole (see Figure \ref{Ts}).
Observe that the graphs in $W_1 \cup W_2 \cup W_3 \cup W_4$ in Figure \ref{Ws} are
tadpoles. A graph in $W_5$ is a tadpole only when it is in $W_3$ and a graph in
$W_6$ is a tadpoles only when it is in $W_2$.

\begin{lemma} \label{s2}
Let $\widehat{G}$ be a connected chordal signed separable bigraph. 
Suppose that no two signed simplicial edges of $\widehat{G}$ induce a $2K_2$. 
Then, for any signed simplicial edge $e$ of $\widehat{G}$, there is an induced
tadpole whose end is incident with $e$.
\end{lemma}
\pf We prove the lemma by induction on the number of vertices of $\widehat{G}$. 
Note that $\widehat{G}$ contains at least six vertices and one can verify easily
that the lemma is true when $\widehat{G}$ has exactly six vertices. So assume 
that $\widehat{G}$ has at least seven vertices and that the lemma is true for all 
connected chordal signed separable bigraphs with fewer vertices than $\widehat{G}$.

Let $S$ be a set of vertices which minimally separates $\widehat{H}_1$ and 
$\widehat{H}_2$. By Lemma \ref{h1orh2}, $e$ is contained in $\widehat{H}_1$ or in 
$\widehat{H}_2$. We assume without loss of generality that $e$ is in $\widehat{H_1}$.
Let $\widehat{G_1}$ (respectively, $\widehat{G_2}$) be the subgraph of $\widehat{G}$ 
induced by $V(\widehat{H_1}) \cup S$ (respectively, $V(\widehat{H_2}) \cup S$).  
Note that every signed simplicial edge of $\widehat{G_2}$ must have an endvertex in 
$S$, as otherwise it is a signed simplicial edge of $\widehat{G}$ which forms 
an induced $2K_2$ with $e$, contradicting the assumption.
Since $\widehat{H_1}$ has at least two vertices, $\widehat{G_2}$ has at most 
$|V(\widehat{G})|-2$ vertices. We prove that it is possible to choose such a set $S$ 
so that $\widehat{G_2}$ is non-separable or has at most $|V(\widehat{G})|-3$ vertices.

For the sake of proof we assume that $\widehat{G_2}$ is separable and has 
$|V(\widehat{G})|-2$ vertices. Then $\widehat{H_1}$ consists of $e$ only and $e$ is 
the only signed simplicial edge of $\widehat{G}$. Since $\widehat{G_2}$ is separable, 
it contains two edges $e'_1, e'_2$ forming an induced $2K_2$. Let $S'$ be
a set which minimally separates $\widehat{H'_1}$ and $\widehat{H'_2}$ where 
$e'_1$ is in $\widehat{H'_1}$ and $e'_2$ is in $\widehat{H'_2}$. 
By Lemma \ref{h1orh2} we can assume that $\widehat{H'_1}$ contains $e$.
Thus $\widehat{H'_1}$ contains both $e$ and $e'_1$, which means it has at least
four vertices and hence the subgraph of $\widehat{G}$ induced by
$V(\widehat{H'_2}) \cup S'$ has at most $|V(\widehat{G})|-4$ vertices.
Therefore we can assume that $\widehat{G_2}$ is non-separable or has at most 
$|V(\widehat{G})|-3$ vertices. 

Suppose that $\widehat{G_2}$ is non-separable. Then by Lemma \ref{s1} $\widehat{G_2}$
contains a graph $W \in W_1 \cup W_2 \cup \cdots \cup W_6$ (see Figure \ref{Ws}) as 
an induced subgraph with the vertices $x, x'$ being the only vertices in $S$.
We claim that $W \notin W_5 \cup W_6$. Indeed, suppose to the contrary that 
$W \in W_5 \cup W_6$ with $x, x'$ being the only vertices in $S$. Then $x, x'$ are 
the only vertices in $W$ and by Lemma \ref{min-ss} they have a common neighbour $y$ 
in $\widehat{H}_1$. Hence the subgraph of $\widehat{G}$ induced by $V(W) \cup \{y\}$
is a graph in $D$ if $W \in W_5$, and a graph in $F_5$ if $W \in W_6$, contradicting
the assumption that $\widehat{G}$ is chordal. Therefore 
$W \in W_1 \cup W_2 \cup W_3 \cup W_4$. If $W \in W_1 \cup W_2$, then 
$W$ together with a shortest path connecting $x$ and $e$ in $\widehat{G}_1$ 
induce a tadpole of type 1. If $W \in W_3\cup W_4$, then 
$W$ together with a shortest path connecting $x$ and $e$ in $\widehat{G}_1$ 
induce a tadpole of type 2. In either case the end of the tadpole is incident with 
the signed simplicial edge $e$.

Suppose now that $\widehat{G_2}$ is separable. Our assumption above ensures that it
has at most $|V(\widehat{G})|-3$ vertices. Since $S$ minimally separates 
$\widehat{H_1}$ and $\widehat{H_2}$, every vertex of $S$ has a neighbour
in $\widehat{H_1}$ and a neighbour in $\widehat{H_2}$. Consider two vertices $s, s'$ 
of $S$ from the same partite set. By Lemma \ref{min-ss}, $s, s'$ have 
a common neighbour in $\widehat{H_2}$ and hence their neighbourhood in $\widehat{H_1}$
are comparable, as otherwise $\widehat{G}$ contains an induced cycle of length 
$\geq 6$, a contradiction to the assumption that $\widehat{G}$ is chordal. 
This implies that the vertices of $S$ of the same partite set have a common neighbour
in $\widehat{H_1}$. Let $u, v$ be a pair of vertices in $\widehat{H_1}$ where $u$ is 
a common neighbour of the vertices of $S$ in one partite set and $v$ is a common
neighbour of the vertices of $S$ in the other partite set. 
(If all vertices of $S$ are in one partite set then let $v$ be any vertex in
$\widehat{H_1}$ adjacent to $u$). 
Let $\widehat{G^*}$ be the graph obtained from the subgraph of $\widehat{G}$ 
induced by $V(\widehat{G_2}) \cup \{u,v\}$ by adding an edge $uv$ (either positive
or negative) if $uv$ is not an edge of $\widehat{G}$. 
Since the neighbourhood of $uv$ is $S$ which forms a positive biclique, $uv$ is a 
signed simplicial edge of $\widehat{G^*}$. The graph graph $\widehat{G^*}$ is chordal.
Indeed, if $uv$ is an edge of $\widehat{G}$, then $\widehat{G^*}$ is an induced 
subgraph of $\widehat{G}$; if $uv$ is not an edge of $\widehat{G}$ then 
$\widehat{G^*}-uv$ is an induced subgraph of $\widehat{G}$.    
Since $\widehat{G_2}$ has at most $|V(\widehat{G})|-3$ vertices, $\widehat{G^*}$ has 
fewer vertices than $\widehat{G}$. Since no two signed simplicial edges of
$\widehat{G}$ induce a $2K_2$, no two signed simplicial edges of $\widehat{G^*}$ 
induce a $2K_2$. Hence, by the inductive hypothesis, $\widehat{G^*}$ contains an 
induced tadpole $T$ with its end incident with $uv$. It is easy to verify that 
$T-\{u,v\}$ is an induced tadpole in $\widehat{G_2}$ with its end being the only
vertex in $S$. The shortest path between the end of $T-\{u,v\}$ and $e$ together with
$T-\{u,v\}$ is an induced tadpole of $\widehat{G}$ with its end incident with $e$.
This completes the proof.
\qed

\begin{cor} \label{s2.5}
Let $\widehat{G}$ be a signed separable bigraph and $S$ be a set which minimally 
separates $\widehat{H_1}$ and $\widehat{H_2}$. Let $\widehat{G_1}$ 
(respectively, $\widehat{G_2}$) be the subgraph of $\widehat{G}$ induced by 
$V(\widehat{H_1}) \cup S$ (respectively, $V(\widehat{H_2}) \cup S$). If either 
of $\widehat{G_1}$, $\widehat{G_2}$ contains a graph in $W_5 \cup W_6$ as an
induced subgraph with $x, x'$ being the only vertices in $S$ then $\widehat{G}$ 
contains a graph in $F_5 \cup D$ as an induced subgraph.
\qed
\end{cor}

\begin{center}
\begin{figure}[htb]
                \center
                \begin{tikzpicture}[>=latex]
                \begin{pgfonlayer}{nodelayer}
\node [style=blackvertex] (1) at (-1,0) {};
\node [style=blackvertex] (2) at (1,0) {};
\node [style=blackvertex] (3) at (2,0) {};
\node [style=blackvertex] (4) at (0,1) {};
\node [style=blackvertex] (5) at (0,-1) {};
\node [style=blackvertex] (6) at (3,0) {};
\node [style=textbox] at (4,0) {$\dots$};
\node [label={below:$x=x'$}] [style=blackvertex] (7) at (5,0) {};

\draw[style=edge] (5) to (2) to (4)[red];
\draw[style=edge] (5) to (1) to (4);
\draw[style=edge] (2) to (3.5,0);
\draw[style=edge] (7) to (4.5,0);

\node [style=blackvertex] (9) at (7,0) {};
\node [style=blackvertex] (10) at (8,0) {};
\node [style=blackvertex] (11) at (9,0) {};
\node [style=blackvertex] (12) at (11,0) {};
\node [style=blackvertex] (13) at (10,1) {};
\node [style=blackvertex] (14) at (10,-1) {};
\node [style=textbox] at (6,0) {$\dots$};

\draw[style=edge] (13) to (11) to (14)[red];
\draw[style=edge] (13) to (12) to (14);
\draw[style=edge] (11) to (6.5,0);
\draw[style=edge] (7) to (5.5,0);

%% (I)

\node [style=blackvertex] (1) at (-1,-3) {};
\node [style=blackvertex] (2) at (1,-3) {};
\node [style=blackvertex] (3) at (2,-3) {};
\node [style=blackvertex] (4) at (0,-2) {};
\node [style=blackvertex] (5) at (0,-4) {};
\node [style=blackvertex] (6) at (3,-3) {};
\node [style=textbox] at (4,-3) {$\dots$};
\node [label={below:$x=x'$}] [style=blackvertex] (7) at (5,-3) {};

\draw[style=edge] (5) to (2) to (4)[red];
\draw[style=edge] (5) to (1) to (4);
\draw[style=edge] (2) to (3.5,-3);
\draw[style=edge] (7) to (4.5,-3);
\draw[style=edge] (7) to (5.5,-3);

\node [style=blackvertex] (9) at (7,-3) {};
\node [style=blackvertex] (10) at (8,-3) {};
\node [style=blackvertex] (11) at (9,-3) {};
\node [style=blackvertex] (12) at (11,-3) {};
\node [style=blackvertex] (13) at (10,-2) {};
\node [style=blackvertex] (14) at (10,-4) {};
\node [style=textbox] at (6,-3) {$\dots$};

\draw[style=edge] (13) to (11) to (14)[red];
\draw[style=edge] (13) to (12) to (14);
\draw[style=edge] (11) to (6.5,-3);
\draw[style=edge] (10) to [out=90,in=180] (9.5,-1.5) to [out=0,in=90] (12)[red];

%%(II)

\node [style=blackvertex] (1) at (-1,-6) {};
\node [style=blackvertex] (2) at (1,-6) {};
\node [style=blackvertex] (3) at (2,-6) {};
\node [style=blackvertex] (4) at (0,-5) {};
\node [style=blackvertex] (5) at (0,-7) {};
\node [style=blackvertex] (6) at (3,-6) {};
\node [style=textbox] at (4,-6) {$\dots$};
\node [label={below:$x=x'$}] [style=blackvertex] (7) at (5,-6) {};

\draw[style=edge] (5) to (2) to (4)[red];
\draw[style=edge] (5) to (1) to (4);
\draw[style=edge] (2) to (3.5,-6);
\draw[style=edge] (7) to (4.5,-6);
\draw[style=edge] (1) to [out=90,in=180] (.5,-4.5) to [out=0,in=90] (3)[red];

\node [style=blackvertex] (9) at (7,-6) {};
\node [style=blackvertex] (10) at (8,-6) {};
\node [style=blackvertex] (11) at (9,-6) {};
\node [style=blackvertex] (12) at (11,-6) {};
\node [style=blackvertex] (13) at (10,-5) {};
\node [style=blackvertex] (14) at (10,-7) {};
\node [style=textbox] at (6,-6) {$\dots$};

\draw[style=edge] (13) to (11) to (14)[red];
\draw[style=edge] (13) to (12) to (14);
\draw[style=edge] (11) to (6.5,-6);
\draw[style=edge] (7) to (5.5,-6);
\draw[style=edge] (10) to [out=90,in=180] (9.5,-4.5) to [out=0,in=90] (12)[red];

%%(III)

                \end{pgfonlayer}
                \end{tikzpicture}
\caption{\label{Sums} Sums of two tadpoles}
        \end{figure}
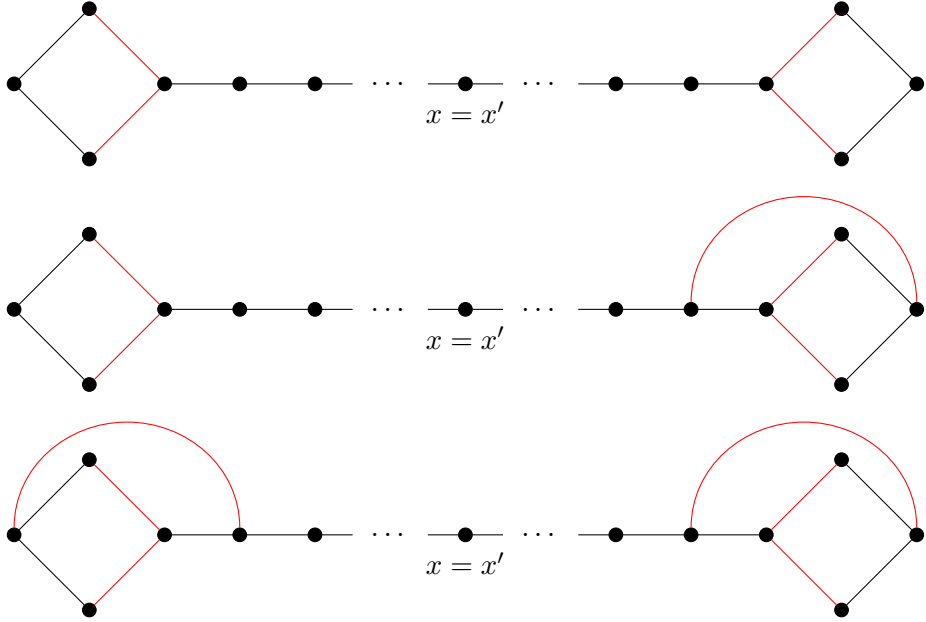
\end{center}

For a tadpole $T$ and a vertex $u$ not in $T$ in a signed bigraph $\widehat{G}$, 
we say that $u$ and $T$ are {\em completely adjacent} (or $u$ is 
{\em completely adjacent to} $T$) if $u$ is adjacent to all vertices of $T$ in 
the partite set of $\widehat{G}$ opposite to $u$.

Let $T$ and $T'$ be two tadpoles with ends $x$ and $x'$ respectively.
The {\em sum} of $T$ and $T'$ is the signed bigraph obtained from $T$ and $T'$ by
identifying $x$ and $x'$. Figure \ref{Sums} depicts all possible sums of tadpoles.
The {\em join} of $T$ and $T'$ is the signed bigraph obtained from $T$ and $T'$ by
\begin{itemize}
\item adding a complete adjacency between $x$ and $T'$ with positive edges, except
      when $T'$ is $W_1$ in which case $xx'$ is an added edge and the other edge 
      between $x$ and the vertex at distance 2 from $x'$ in $T'$ may be either 
      positive or negative, and

\item adding a complete adjacency between $x'$ and $T$ with positive edges, except
      when $T$ is $W_1$ in which case $xx'$ is an added edge and the other edge 
      between $x'$ and the vertex at distance 2 from $x$ in $T$ may be either 
      positive or negative. 
\end{itemize}

\begin{center}
\begin{figure}[htb!]
                \center
                \begin{tikzpicture}[>=latex]
                \begin{pgfonlayer}{nodelayer}
\node [style=blackvertex] (1) at (-1,0) {};
\node [style=blackvertex] (2) at (0,1) {};
\node [style=blackvertex] (3) at (0,-1) {};
\node [style=blackvertex] (4) at (1,0) {};
\node [style=blackvertex] (5) at (2,0) {};
\node [style=blackvertex] (6) at (3,0) {};
\node [style=textbox] at (4,0) {$\dots$};
\node [style=blackvertex] (8) at (5,0) {};
\node [style=blackvertex] (9) at (6,0) {};
\node [label={above:$x$}] [style=blackvertex] (10) at (7,1) {};

\draw[style=edge] (2) to (4) to (3)[red];
\draw[style=edge] (2) to (1) to (3);
\draw[style=edge] (4) to (3.5,0);
\draw[style=edge] (10) to (9) to (4.5,0);

\node [label={below:$x'$}] [style=blackvertex] (11) at (7,-1) {};
\node [style=blackvertex] (12) at (8,0) {};
\node [style=blackvertex] (13) at (9,0) {};
\node [style=textbox] at (10,0) {$\dots$};
\node [style=blackvertex] (15) at (11,0) {};
\node [style=blackvertex] (16) at (12,0) {};
\node [style=blackvertex] (17) at (13,0) {};
\node [style=blackvertex] (18) at (14,1) {};
\node [style=blackvertex] (19) at (14,-1) {};
\node [style=blackvertex] (20) at (15,0) {};

\draw[style=edge] (18) to (17) to (19)[red];
\draw[style=edge] (18) to (20) to (19);
\draw[style=edge] (17) to (10.5,0);
\draw[style=edge] (11) to (12) to (9.5,0);

\draw[style=edge] (16) to (10) to (12)[blue];
\draw[style=edge] (9.5,.2) to (10) to (10,.2)[blue];
\draw[style=edge] (10) to [out=5,in=170] (18)[blue];
\draw[style=edge] (10) to [out=-8,in=110] (19)[blue];

\draw[style=edge] (5) to (11) to (9)[blue];
\draw[style=edge] (4,-.2) to (11) to (4.5,-.2)[blue];
\draw[style=edge] (11) to [out=185,in=-10] (3)[blue];
\draw[style=edge] (11) to [out=172,in=-70] (2)[blue];

%% (I)

\node [style=blackvertex] (1) at (-1,-3.5) {};
\node [style=blackvertex] (2) at (0,-2.5) {};
\node [style=blackvertex] (3) at (0,-4.5) {};
\node [style=blackvertex] (4) at (1,-3.5) {};
\node [style=blackvertex] (5) at (2,-3.5) {};
\node [style=blackvertex] (6) at (3,-3.5) {};
\node [style=textbox] at (4,-3.5) {$\dots$};
\node [style=blackvertex] (8) at (5,-3.5) {};
\node [style=blackvertex] (9) at (6,-3.5) {};
\node [label={above:$x$}] [style=blackvertex] (10) at (7,-2.5) {};

\draw[style=edge] (2) to (4) to (3)[red];
\draw[style=edge] (2) to (1) to (3);
\draw[style=edge] (4) to (3.5,-3.5);
\draw[style=edge] (10) to (9) to (4.5,-3.5);

\node [label={below:$x'$}] [style=blackvertex] (11) at (7,-4.5) {};
\node [style=blackvertex] (12) at (8,-3.5) {};
\node [style=blackvertex] (13) at (9,-3.5) {};
\node [style=textbox] at (10,-3.5) {$\dots$};
\node [style=blackvertex] (15) at (11,-3.5) {};
\node [style=blackvertex] (16) at (12,-3.5) {};
\node [style=blackvertex] (17) at (13,-3.5) {};
\node [style=blackvertex] (18) at (14,-2.5) {};
\node [style=blackvertex] (19) at (14,-4.5) {};
\node [style=blackvertex] (20) at (15,-3.5) {};

\draw[style=edge] (18) to (17) to (19)[red];
\draw[style=edge] (18) to (20) to (19);
\draw[style=edge] (17) to (10.5,-3.5);
\draw[style=edge] (11) to (12) to (9.5,-3.5);

\draw[style=edge] (17) to (10) to (15)[blue];
\draw[style=edge] (12) to (10) to (9.5,-3.3)[blue];
\draw[style=edge] (10) to [out=10,in=90] (20)[blue];

\draw[style=edge] (5) to (11) to (9)[blue];
\draw[style=edge] (4,-3.7) to (11) to (4.5,-3.7)[blue];
\draw[style=edge] (11) to [out=185,in=-10] (3)[blue];
\draw[style=edge] (11) to [out=172,in=-70] (2)[blue];

%%(II)

\node [style=blackvertex] (1) at (-1,-7) {};
\node [style=blackvertex] (2) at (0,-6) {};
\node [style=blackvertex] (3) at (0,-8) {};
\node [style=blackvertex] (4) at (1,-7) {};
\node [style=blackvertex] (5) at (2,-7) {};
\node [style=blackvertex] (6) at (3,-7) {};
\node [style=textbox] at (4,-7) {$\dots$};
\node [style=blackvertex] (8) at (5,-7) {};
\node [style=blackvertex] (9) at (6,-7) {};
\node [label={above:$x$}] [style=blackvertex] (10) at (7,-6) {};

\draw[style=edge] (2) to (4) to (3)[red];
\draw[style=edge] (2) to (1) to (3);
\draw[style=edge] (4) to (3.5,-7);
\draw[style=edge] (10) to (9) to (4.5,-7);

\node [label={below:$x'$}] [style=blackvertex] (11) at (7,-8) {};
\node [style=blackvertex] (12) at (8,-7) {};
\node [style=blackvertex] (13) at (9,-7) {};
\node [style=textbox] at (10,-7) {$\dots$};
\node [style=blackvertex] (15) at (11,-7) {};
\node [style=blackvertex] (16) at (12,-7) {};
\node [style=blackvertex] (17) at (13,-7) {};
\node [style=blackvertex] (18) at (14,-6) {};
\node [style=blackvertex] (19) at (14,-8) {};
\node [style=blackvertex] (20) at (15,-7) {};

\draw[style=edge] (18) to (17) to (19)[red];
\draw[style=edge] (18) to (20) to (19);
\draw[style=edge] (17) to (10.5,-7);
\draw[style=edge] (11) to (12) to (9.5,-7);

\draw[style=edge] (17) to (10) to (15)[blue];
\draw[style=edge] (12) to (10) to (9.5,-6.8)[blue];
\draw[style=edge] (10) to [out=10,in=90] (20)[blue];

\draw[style=edge] (4) to (11) to (6) [blue];
\draw[style=edge] (9) to (11) to (4.5,-7.2)[blue];
\draw[style=edge] (11) to [out=190,in=-90] (1)[blue];

%%(III)

                \end{pgfonlayer}
                \end{tikzpicture}
\caption{\label{Joins1} Joins of tadpoles $T, T'$ of type 1 (with ends $x, x'$ 
respectively) where $xx'$ is not an edge, $x$ is completely adjacent to $T'$ by 
positive edges, and $x'$ is completely adjacent to $T$ by positive edges}
        \end{figure}
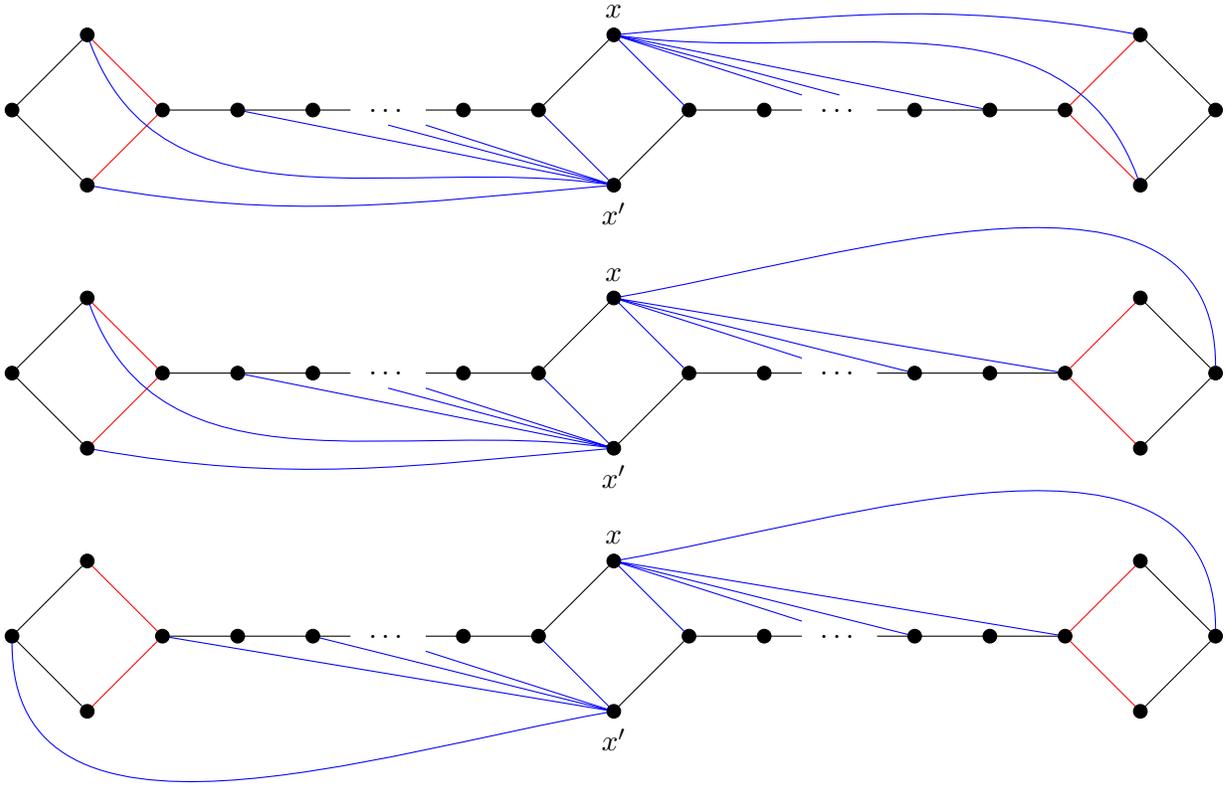
\end{center}

For the space reason, we only show in drawing the joins of tapoles of type 1.
Figure \ref{Joins1} depicts the joins of tadpoles $T, T'$ of type 1 where $xx'$ is 
not an edge and there is a complete adjacency with positive edges between $x$ and 
$T'$ and between $x'$ and $T$. Figure \ref{Joins2} depicts the joins of tadpoles 
$T, T'$ of type 1 where $xx'$ is an edge and there is a complete adjacency with 
positive edges between $x$ and $T'$ and between $x'$ and $T$.
Figure \ref{JoinsW1} depicts the joins of $T \in W_1$ and $T'$ of type 1 (which may
also be in $W_1$) where the edge between $x'$ and the vertex at distance 2 from $x$ 
in $T$ is negative (in the case when $T, T' \in  W_1$ the edge between $x$ and 
the vertex at distance 2 from $x'$ in $T'$ may also be negative as shown in the last
one).

\begin{center}
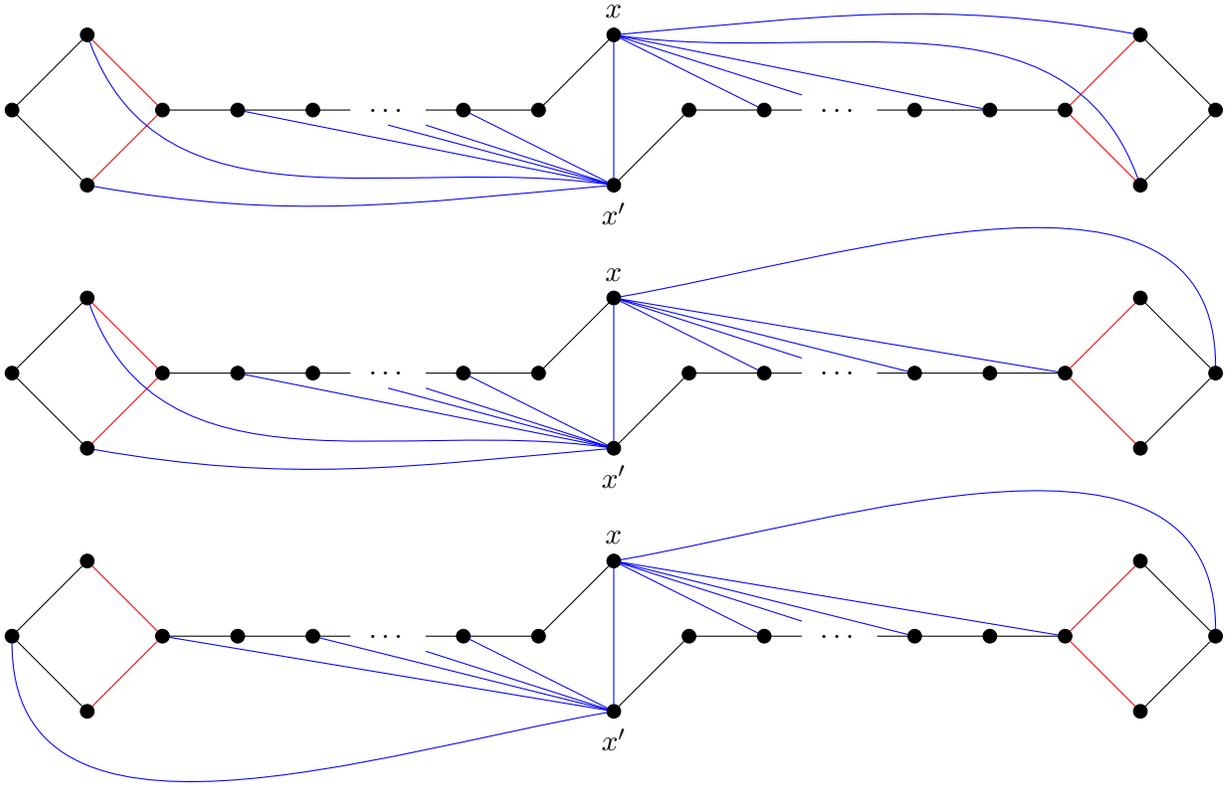
\begin{figure}[htb!]
                \center
                \begin{tikzpicture}[>=latex]
                \begin{pgfonlayer}{nodelayer}

\node [style=blackvertex] (1) at (-1,-10.5) {};
\node [style=blackvertex] (2) at (0,-9.5) {};
\node [style=blackvertex] (3) at (0,-11.5) {};
\node [style=blackvertex] (4) at (1,-10.5) {};
\node [style=blackvertex] (5) at (2,-10.5) {};
\node [style=blackvertex] (6) at (3,-10.5) {};
\node [style=textbox] at (4,-10.5) {$\dots$};
\node [style=blackvertex] (8) at (5,-10.5) {};
\node [style=blackvertex] (9) at (6,-10.5) {};
\node [label={above:$x$}] [style=blackvertex] (10) at (7,-9.5) {};

\draw[style=edge] (2) to (4) to (3)[red];
\draw[style=edge] (2) to (1) to (3);
\draw[style=edge] (4) to (3.5,-10.5);
\draw[style=edge] (10) to (9) to (4.5,-10.5);

\node [label={below:$x'$}] [style=blackvertex] (11) at (7,-11.5) {};
\node [style=blackvertex] (12) at (8,-10.5) {};
\node [style=blackvertex] (13) at (9,-10.5) {};
\node [style=textbox] at (10,-10.5) {$\dots$};
\node [style=blackvertex] (15) at (11,-10.5) {};
\node [style=blackvertex] (16) at (12,-10.5) {};
\node [style=blackvertex] (17) at (13,-10.5) {};
\node [style=blackvertex] (18) at (14,-9.5) {};
\node [style=blackvertex] (19) at (14,-11.5) {};
\node [style=blackvertex] (20) at (15,-10.5) {};

\draw[style=edge] (18) to (17) to (19)[red];
\draw[style=edge] (18) to (20) to (19);
\draw[style=edge] (17) to (10.5,-10.5);
\draw[style=edge] (11) to (12) to (9.5,-10.5);

\draw[style=edge] (11) to (10) to (13)[blue];
\draw[style=edge] (16) to (10) to (9.5,-10.3)[blue];
\draw[style=edge] (10) to [out=5,in=170] (18)[blue];
\draw[style=edge] (10) to [out=-8,in=110] (19)[blue];

\draw[style=edge] (8) to (11) to (5) [blue];
\draw[style=edge] (4,-10.7) to (11) to (4.5,-10.7)[blue];
\draw[style=edge] (11) to [out=185,in=-10] (3)[blue];
\draw[style=edge] (11) to [out=172,in=-70] (2)[blue];

%%(IV)

\node [style=blackvertex] (1) at (-1,-14) {};
\node [style=blackvertex] (2) at (0,-13) {};
\node [style=blackvertex] (3) at (0,-15) {};
\node [style=blackvertex] (4) at (1,-14) {};
\node [style=blackvertex] (5) at (2,-14) {};
\node [style=blackvertex] (6) at (3,-14) {};
\node [style=textbox] at (4,-14) {$\dots$};
\node [style=blackvertex] (8) at (5,-14) {};
\node [style=blackvertex] (9) at (6,-14) {};
\node [label={above:$x$}] [style=blackvertex] (10) at (7,-13) {};

\draw[style=edge] (2) to (4) to (3)[red];
\draw[style=edge] (2) to (1) to (3);
\draw[style=edge] (4) to (3.5,-14);
\draw[style=edge] (10) to (9) to (4.5,-14);

\node [label={below:$x'$}] [style=blackvertex] (11) at (7,-15) {};
\node [style=blackvertex] (12) at (8,-14) {};
\node [style=blackvertex] (13) at (9,-14) {};
\node [style=textbox] at (10,-14) {$\dots$};
\node [style=blackvertex] (15) at (11,-14) {};
\node [style=blackvertex] (16) at (12,-14) {};
\node [style=blackvertex] (17) at (13,-14) {};
\node [style=blackvertex] (18) at (14,-13) {};
\node [style=blackvertex] (19) at (14,-15) {};
\node [style=blackvertex] (20) at (15,-14) {};

\draw[style=edge] (18) to (17) to (19)[red];
\draw[style=edge] (18) to (20) to (19);
\draw[style=edge] (17) to (10.5,-14);
\draw[style=edge] (11) to (12) to (9.5,-14);

\draw[style=edge] (11) to (10) to (13)[blue];
\draw[style=edge] (15) to (10) to (17)[blue];
\draw[style=edge] (10) to (9.5,-13.8)[blue];
\draw[style=edge] (10) to [out=10,in=90] (20)[blue];

\draw[style=edge] (8) to (11) to (5) [blue];
\draw[style=edge] (4,-14.2) to (11) to (4.5,-14.2)[blue];
\draw[style=edge] (11) to [out=185,in=-10] (3)[blue];
\draw[style=edge] (11) to [out=172,in=-70] (2)[blue];

%%(V)

\node [style=blackvertex] (1) at (-1,-17.5) {};
\node [style=blackvertex] (2) at (0,-16.5) {};
\node [style=blackvertex] (3) at (0,-18.5) {};
\node [style=blackvertex] (4) at (1,-17.5) {};
\node [style=blackvertex] (5) at (2,-17.5) {};
\node [style=blackvertex] (6) at (3,-17.5) {};
\node [style=textbox] at (4,-17.5) {$\dots$};
\node [style=blackvertex] (8) at (5,-17.5) {};
\node [style=blackvertex] (9) at (6,-17.5) {};
\node [label={above:$x$}] [style=blackvertex] (10) at (7,-16.5) {};

\draw[style=edge] (2) to (4) to (3)[red];
\draw[style=edge] (2) to (1) to (3);
\draw[style=edge] (4) to (3.5,-17.5);
\draw[style=edge] (10) to (9) to (4.5,-17.5);

\node [label={below:$x'$}] [style=blackvertex] (11) at (7,-18.5) {};
\node [style=blackvertex] (12) at (8,-17.5) {};
\node [style=blackvertex] (13) at (9,-17.5) {};
\node [style=textbox] at (10,-17.5) {$\dots$};
\node [style=blackvertex] (15) at (11,-17.5) {};
\node [style=blackvertex] (16) at (12,-17.5) {};
\node [style=blackvertex] (17) at (13,-17.5) {};
\node [style=blackvertex] (18) at (14,-16.5) {};
\node [style=blackvertex] (19) at (14,-18.5) {};
\node [style=blackvertex] (20) at (15,-17.5) {};

\draw[style=edge] (18) to (17) to (19)[red];
\draw[style=edge] (18) to (20) to (19);
\draw[style=edge] (17) to (10.5,-17.5);
\draw[style=edge] (11) to (12) to (9.5,-17.5);

\draw[style=edge] (11) to (10) to (13)[blue];
\draw[style=edge] (15) to (10) to (17)[blue];
\draw[style=edge] (10) to (9.5,-17.3)[blue];
\draw[style=edge] (10) to [out=10,in=90] (20)[blue];

\draw[style=edge] (8) to (11) to (6) [blue];
\draw[style=edge] (4) to (11) to (4.5,-17.7)[blue];
\draw[style=edge] (11) to [out=190,in=-90] (1)[blue];

%%(VI)

                \end{pgfonlayer}
                \end{tikzpicture}
\caption{\label{Joins2} Joins of tadpoles $T, T'$ of type 1 (with ends $x, x'$ 
respectively) where $xx'$ is an edge, $x$ is completely adjacent to $T'$ by positive 
edges, and $x'$ is completely adjacent to $T$ by positive edges}
        \end{figure}
\end{center}

\begin{center}
\begin{figure}[htb!]
                \center
                \begin{tikzpicture}[>=latex]
                \begin{pgfonlayer}{nodelayer}
\node [style=blackvertex] (1) at (-1,1) {};
\node [style=blackvertex] (2) at (0,2) {};
\node [style=blackvertex] (3) at (0,0) {};
\node [label={above:$x$}] [style=blackvertex] (4) at (1,1) {};
\node [label={below:$x'$}] [style=blackvertex] (5) at (2,0) {};
\node [style=blackvertex] (6) at (3,0) {};
\node [style=blackvertex] (7) at (4,0) {};
\node [style=blackvertex] (8) at (6,0) {};
\node [style=blackvertex] (9) at (7,0) {};
\node [style=blackvertex] (10) at (8,1) {};
\node [style=blackvertex] (11) at (8,-1) {};
\node [style=blackvertex] (12) at (9,0) {};
\node [style=textbox] at (5,0) {$\dots$};

\draw[style=edge] (2) to (1) to (3);
\draw[style=edge] (2) to (4) to (3)[red];
\draw[style=edge] (5) to (6) to (7) to (4.5,0);
\draw[style=edge] (5.5,0) to (8) to (9);
\draw[style=edge] (10) to (12) to (11);
\draw[style=edge] (10) to (9) to (11)[red];

\draw[style=edge] (1) to [out=-90,in=240] (5)[red];
\draw[style=edge] (8) to (4) to (4.5,.1)[blue];
\draw[style=edge] (7) to (4) to (5)[blue];
\draw[style=edge] (4) to [out=0,in=90] (11)[blue];
\draw[style=edge] (4) to [out=45,in=120] (10)[blue];

%(I)

\node [style=blackvertex] (1) at (-1,-3) {};
\node [style=blackvertex] (2) at (0,-2) {};
\node [style=blackvertex] (3) at (0,-4) {};
\node [label={above:$x$}] [style=blackvertex] (4) at (1,-3) {};
\node [label={below:$x'$}] [style=blackvertex] (5) at (2,-4) {};
\node [style=blackvertex] (6) at (3,-4) {};
\node [style=blackvertex] (7) at (4,-4) {};
\node [style=blackvertex] (8) at (6,-4) {};
\node [style=blackvertex] (9) at (7,-4) {};
\node [style=blackvertex] (10) at (8,-3) {};
\node [style=blackvertex] (11) at (8,-5) {};
\node [style=blackvertex] (12) at (9,-4) {};
\node [style=textbox] at (5,-4) {$\dots$};

\draw[style=edge] (2) to (1) to (3);
\draw[style=edge] (2) to (4) to (3)[red];
\draw[style=edge] (5) to (6) to (7) to (4.5,-4);
\draw[style=edge] (5.5,-4) to (8) to (9);
\draw[style=edge] (10) to (12) to (11);
\draw[style=edge] (10) to (9) to (11)[red];

\draw[style=edge] (1) to [out=-90,in=240] (5)[red];
\draw[style=edge] (7) to (4) to (5)[blue];
\draw[style=edge] (9) to (4) to (4.5,-3.9)[blue];
\draw[style=edge] (4) to [out=45,in=100] (12)[blue];

%(II)

\node [style=blackvertex] (1) at (-2,-7) {};
\node [style=blackvertex] (2) at (-1,-6) {};
\node [style=blackvertex] (3) at (-1,-8) {};
\node [label={above:$x$}] [style=blackvertex] (4) at (0,-7) {};
\node [label={below:$x'$}] [style=blackvertex] (5) at (1,-8) {};
\node [style=blackvertex] (6) at (2,-7) {};
\node [style=blackvertex] (7) at (2,-9) {};
\node [style=blackvertex] (8) at (3,-8) {};

\draw[style=edge] (2) to (1) to (3);
\draw[style=edge] (2) to (4) to (3)[red];
\draw[style=edge] (6) to (8) to (7);
\draw[style=edge] (6) to (5) to (7)[red];

\draw[style=edge] (1) to [out=-90,in=240] (5)[red];
\draw[style=edge] (4) to (5)[blue];
\draw[style=edge] (4) to [out=45,in=90] (8)[blue];

%(III)

\node [style=blackvertex] (1) at (5,-7) {}; 
\node [style=blackvertex] (2) at (6,-6) {};
\node [style=blackvertex] (3) at (6,-8) {};
\node [label={above:$x$}] [style=blackvertex] (4) at (7,-7) {};
\node [label={below:$x'$}] [style=blackvertex] (5) at (8,-8) {};
\node [style=blackvertex] (6) at (9,-7) {};
\node [style=blackvertex] (7) at (9,-9) {};
\node [style=blackvertex] (8) at (10,-8) {};

\draw[style=edge] (2) to (1) to (3);
\draw[style=edge] (2) to (4) to (3)[red];
\draw[style=edge] (6) to (8) to (7);
\draw[style=edge] (6) to (5) to (7)[red];

\draw[style=edge] (1) to [out=-90,in=240] (5)[red];
\draw[style=edge] (4) to (5)[blue];
\draw[style=edge] (4) to [out=45,in=90] (8)[red];

%(IV)

                \end{pgfonlayer}
                \end{tikzpicture}
\caption{\label{JoinsW1} Joins of a graph in $W_1$ (with end $x$) and tadpole $T'$ of
type 1 
(with end $x'$) where the edge between $x'$ and the vertex at distance 2 from $x$ in
$W_1$ is negative; when $T' \in W_1$, the edge between $x$ and the vertex 
at distance 2 from $x'$ in $T'$ is negative (see the last graph)}
        \end{figure}
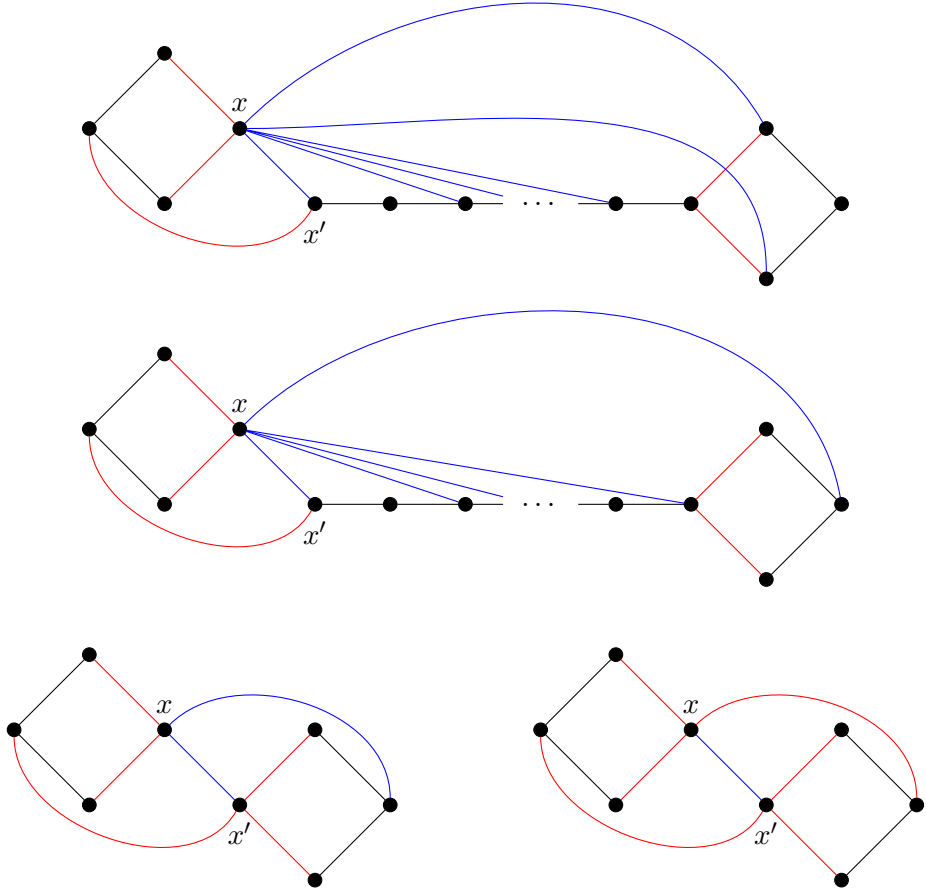
\end{center}

\begin{center}
\begin{figure}[htb]
                \center
                \begin{tikzpicture}[>=latex]
                \begin{pgfonlayer}{nodelayer}
\node [label={left:$w$}] [style=blackvertex] (1) at (-1,0) {};
\node [label={below:$x_k$}] [style=blackvertex] (2) at (1,0) {};
\node [label={below:$x_{k-1}$}] [style=blackvertex] (3) at (2,0) {};
\node [label={above:$y$}] [style=blackvertex] (4) at (0,1) {};
\node [label={below:$z$}] [style=blackvertex] (5) at (0,-1) {};
\node [style=blackvertex] (6) at (3,0) {};
\node [style=textbox] at (4,0) {$\dots$};
\node [label={below:$x=x_1$}] [style=blackvertex] (7) at (5,0) {};

\node [style=textbox] at (2.5,-1.5) {${\cal T}_1$};

\draw[style=edge] (5) to (2) to (4)[red];
\draw[style=edge] (5) to (1) to (4);
\draw[style=edge] (2) to (3.5,0);
\draw[style=edge] (7) to (4.5,0);
%\draw[style=edge] (1) to [out=90,in=180] (4.5,1.5) to [out=0,in=90] (3)[red];

%(I)
                
\node [label={left:$w$}] [style=blackvertex] (1) at (7,0) {};
\node [label={below:$x_k$}] [style=blackvertex] (2) at (9,0) {};
\node [label={below:$x_{k-1}$}] [style=blackvertex] (3) at (10,0) {};
\node [label={above:$y$}] [style=blackvertex] (4) at (8,1) {};
\node [label={below:$z$}] [style=blackvertex] (5) at (8,-1) {};
\node [style=blackvertex] (6) at (11,0) {};
\node [style=textbox] at (12,0) {$\dots$};
\node [label={below:$x=x_1$}] [style=blackvertex] (7) at (13,0) {};

\node [style=textbox] at (10.5,-1.5) {${\cal T}_2$};

\draw[style=edge] (5) to (2) to (4)[red];
\draw[style=edge] (5) to (1) to (4);
\draw[style=edge] (2) to (11.5,0); 
\draw[style=edge] (7) to (12.5,0);
\draw[style=edge] (1) to [out=90,in=180] (8.4,1.8) to [out=0,in=90] (3)[red];

%(II)
                \end{pgfonlayer}
                \end{tikzpicture}
\caption{\label{labeled-Ts} Labeled tadpoles}
        \end{figure}
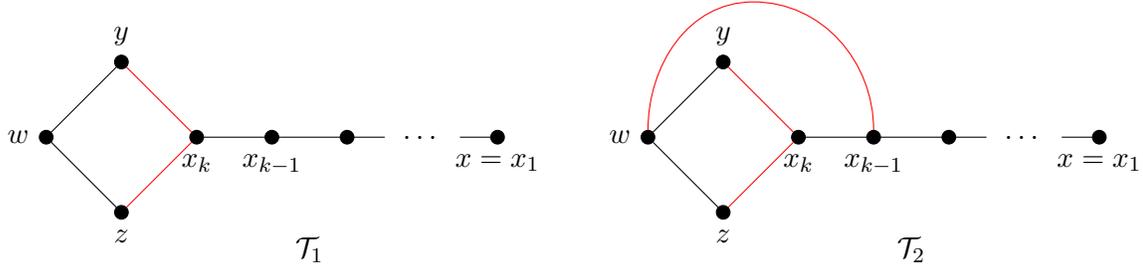
\end{center}

Denote by $\cal S$ (respectively, $\cal J$) the set of sums (respectively, joins)
of tadpoles, and ${\cal F} = F_1 \cup F_2 \cup \cdots \cup F_6 \cup {\cal C} \cup
D \cup {\cal S} \cup {\cal J}$. 
We say that a signed graph is {\em ${\cal F}$-free} if it does not contain any graph
in $\cal F$ as an induced subgraph.

For convenience we label the vertices of tadpoles as shown in
Figure~\ref{labeled-Ts}. Call the vertices $w, y, z$ in a tadpole the {\em heads} of
the tadpole.

\begin{lemma} \label{s4.5}
Let $\widehat{G}$ be an ${\cal F}$-free signed bigraph.
Let $T$ be an induced tadpole in $\widehat{G}$ 
whose vertices are labeled as in Figure~\ref{labeled-Ts} (with $x_1 = x$ being its 
end). Let $x'$ be a vertex not on $T$ such that
\begin{itemize}
\item either $xx'$ is an edge, or
\item $x, x'$ have a common neighbour adjacent to no vertex of $T$ other
      than $x$.
\end{itemize}
If $x'$ is adjacent to a head of $T$ then $x'$ is completely adjacent to $T$ by 
positive edges, except when $x'x$ or $x'w$ is an edge in which case either one 
can be positive or negative.
\end{lemma}
\pf Suppose first that $xx'$ is an edge. Consider first the case when $k$ is
odd. Since $x'$ is adjacent to a head of $T$, $x'w$ is an edge. Thus the statement 
holds when $k = 1$. So assume that $k \geq 3$. Let $x'x_t$ be an edge with $t$ 
being the largest. We claim that $t = k$. Suppose not; $t < k$. When $T$ is of type 1,
$x'x_tx_{t+1} \dots x_kywx'$ is an induced cycle of length $\geq 6$, a contradiction 
to the assumption. When $T'$ is of type 2,  if $t = k-2$, 
then $x', x_{k-2}, x_{k-1}, x_k, y, w$ induce a graph in $D$, otherwise 
$t < k-2$, $x'x_tx_{t+1}\dots x_{k-1}wx'$ is an induced cycle of length $\geq 6$, contradictions. Hence $x'x_k$ is an edge. We must have that $x'$ is completely
adjacent to $T$ as otherwise there is a cycle of length $\geq 6$. 
We show that $x'x_i$ is positive for all $i \geq 2$. Suppose not; $x'x_i$ is negative for 
any $i$ with $i \geq 2$. If $i < k$ then $\widehat{G}$ contains a graph in $D$ 
induced by $x', x_{i-2}, x_{i-1}, x_i, x_{i+1}, x_{i+2}$. If $i = k$ then 
the subgraph of $\widehat{G}$ induced by $x', x_{k-2}, x_{k-1}, x_k, y, w$ 
is a graph in $D$ when $T$ is of type 1, or a graph in $F_5$ when $T$ is of type 2.

Consider then the case when $k$ is even. By assumption and the symmetry between $y$
and $z$ we assume that $x'$ is adjacent to $y$. 
Let $x'x_r$ be an edge with $r$ being the largest. We claim that $r = k-1$ as 
otherwise $x'x_rx_{r+1} \dots x_kyx'$ is an induced cycle of length $\geq 6$, 
a contradiction. The edge $x'z$ must be present in $\widehat{G}$ as otherwise
$x', x_{k-1}, x_k, z, w, y$ induce a graph in $D$ if $T$ is of type 1, or a graph 
in $F_5$ if $T$ is of type 2. In fact, both $x'y, x'z$ must be  
positive or else the same six vertices $x', x_{k-1}, x_k, z, w, y$ induce a 
graph in $F_5$ when $T$ is of type 1, or a graph in $F_4$ when $T$ is of type 2. 
It is now easy to see that $x'$ is completely adjacent to $T$ by positive edges 
except possibly $x'x$.

Suppose now that $u$ is a common neighbour of $x, x'$ which is adjacent to no vertex 
of $T$ other than $x$. As above we consider two cases depending on the parity of $k$. 
Assume first that $k$ is odd. By the symmetry between $y$ and $z$ we assume that $x'$ 
is adjacent to $y$. In the case when $k = 1$, $x'z$ must be present and both 
$x'y, x'z$ are positive, as otherwise $u, x', x, y, z, w$ induce a graph in $D$ or in
$F_5$. So $k \geq 3$. There must be a complete adjacency between $x'$ and $T$ by
positive edges. Indeed, if $x'$ is not completely adjacent to $T$ by positive
edges among the vertices $x_i$ then $\widehat{G}$ contains an induced cycle of length
$\geq 6$ or an induced graph in $D$; if $x'$ is not adjacent to $z$ or either of 
$x'y, x'z$ is negative then $\widehat{G}$ contains an induced graph in $F_5 \cup D$
when $T$ is of type 1 and an induced graph in $F_4 \cup F_5$ when $T$ is of type 2.

Finally, assume that $k$ is even. Then the only head of $T$ adjacent to $x'$ is $w$.
A similar argument as above shows that $x'$ is completely adjacent to $T$ by 
positive edges except possibly $x'w$ or else $\widehat{G}$ contains a graph 
in ${\cal C} \cup D \cup F_5$ as an induced subgraph.
This completes the proof.
\qed

\section{Main theorem} \label{5}

Recall that 
$${\cal F} = F_1 \cup F_2 \cup \cdots \cup F_6 \cup {\cal C} \cup D \cup {\cal S} 
              \cup {\cal J}.$$ 
We make two remarks on the graphs in $\cal F$.
First, since none of the graphs in $\cal F$ has a signed simplicial edge, none of 
them is chordal and hence can not be an induced subgraph of a chordal signed bigraph. 
Second, if a signed bigraph $\widehat{G}$ is $\cal F$-free and $e$ is a signed 
simplicial edge in $\widehat{G}$ then $\widehat{G}-e$ is again $\cal F$-free. 

\begin{theorem} \label{main}
Let $\widehat{G}$ be a signed bigraph. Then $\widehat{G}$ is chordal if and only if
it is ${\cal F}$-free.
\end{theorem}  
\pf The necessity has been discussed above so we only prove the sufficiency.
Suppose to the contrary that there exists an $\cal F$-free signed bigraph which 
is not chordal. Let $\widehat{G}$ be such a graph that is minimal with respect to
vertex deletion, that is, deleting a vertex from $\widehat{G}$ results in a chordal 
signed bigraph. In view of the remarks above we can assume that $\widehat{G}$ 
contains no signed simplicial edge. Theorem \ref{non-separableForbidden}
implies that $\widehat{G}$ is separable.  

Let $S$ be a set in $\widehat{G}$ that minimally separates $\widehat{H_1}$ 
and $\widehat{H_2}$. Let $\widehat{G_1}$ (respectively, $\widehat{G_2}$) be the 
subgraph of $\widehat{G}$ induced by $V(\widehat{H_1}) \cup S$ (respectively,
by $V(\widehat{H_2}) \cup S$). The minimality of $\widehat{G}$ ensures that 
both $\widehat{G_1}$ and $\widehat{G_2}$ are chordal (and hence each has a 
signed simplicial edge). Since $\widehat{G}$ does not have a signed simplicial edge, 
any signed simplicial edge of $\widehat{G_1}$ or $\widehat{G_2}$ must have an
endvertex in $S$. 

We show that $\widehat{G_1}$ and $\widehat{G_2}$ each contains an induced
tadpole with its end being the only vertex in $S$. By symmetry we only prove it 
for $\widehat{G_1}$. If $\widehat{G_1}$ is non-separable, then,  
by Lemma \ref{s1} and Corollary \ref{s2.5}, $\widehat{G_1}$ contains an induced 
tadpole in $W_1 \cup W_2 \cup W_3 \cup W_4$ with its end being the only
vertex in $S$. Suppose that $\widehat{G_1}$ is separable. Since $\widehat{H_2}$ has 
an edge, it has at least two vertices. If $\widehat{H_2}$ has exactly two vertices
then the (only) edge in $\widehat{H_2}$ is a signed simplicial edge of $\widehat{G}$,
a contradiction to the choice of $\widehat{G}$. Thus $\widehat{H_2}$ has at least
three vertices, which means that $\widehat{G_1}$ has at most $|V(\widehat{G})|-3$
vertices. Let $u, v$ be a pair of vertices in $\widehat{H_2}$ where $u$ is a common
neighbour of the vertices of $S$ in one partite set and $v$ is a common neighbour of 
the vertices of $S$ in the other partite set. 
(If all vertices of $S$ are in one partite set then let $v$ be any vertex in
$\widehat{H_1}$ adjacent to $u$).
Let $\widehat{G^*}$ be the graph obtained from the subgraph of $\widehat{G}$
induced by $V(\widehat{G_1}) \cup \{u,v\}$ by adding an edge $uv$ (either positive
or negative) if $uv$ is not an edge of $\widehat{G}$.
Since the neighbourhood of $uv$ is $S$ which forms a positive biclique, $uv$ is a
signed simplicial edge of $\widehat{G^*}$. Either $\widehat{G^*}$ or 
$\widehat{G^*} -uv$ is an induced subgraph of $\widehat{G}$ with fewer vertices
than $\widehat{G}$, $\widehat{G^*}$ is chordal
and $uv$ is the only signed simplicial edge of $\widehat{G^*}$.
By Lemma \ref{s2}, $\widehat{G^*}$
contains an induced tadpole whose end is incident with $uv$. By deleting $u, v$
from the tadpole we obtain a tadpole $T$ in $\widehat{G_1}$ with its end being 
the only vertex in $S$.  Similarly, $\widehat{G_2}$ contains an induced tadpole 
$T'$ with its end being the only vertex in $S$.

Let $x_1, x_2, \dots, x_k, w, y, z$ be the vertices of $T$ as depicted in 
Figure \ref{labeled-Ts}. In a similar way, let 
$x'_1, x'_2, \dots, x'_{\ell}, w', y', z'$ be the vertices of $T'$. Note that 
$x_1$ and $x'_1$ are the ends of $T$ and $T'$ respectively. We know from above that
they are the only vertices in $S$. 

We will show that $\widehat{G}$ contains a graph in $\cal F$ as an induced subgraph 
which contradicts the assumption that $\widehat{G}$ is $\cal F$-free. Note that
any edge between $T$ and $T'$ is incident with at least one of $x_1$ and $x'_1$. 
If $x_1 = x'_1$, then $T$ and $T'$ together induce a graph in $\cal S$, 
contradicting the assumption that $\widehat{G}$ is $\cal F$-free. So $x_1 \neq x'_1$.

If $x'_1$ is not adjacent to any of $w, y, z$ and is adjacent to $x_j$, choose 
such $j$ to be the largest, then $T - \{x_1, x_2, \dots, x_{j-1}\}$ and $T'$ together
induce a graph in $\cal S$, a contradiction. So we may assume that $x'_1$ must 
be adjacent to one of $w, y, z$ if it is adjacent to $x_j$ for some $j$. 
Similarly, $x_1$ must be adjacent to one of $w', y', z'$ if it is adjacent to 
$x'_j$ for some $j$. 

Suppose first that $x_1$ and $x'_1$ are in different partite set of $\widehat{G}$.
Then they are joint by a positive edge as $S$ induces a positive biclique in 
$\widehat{G}$. If $x_1x'_1$ is the only edge between $T$ and $T'$ then again 
$T$ and $T'$ induce a graph in $\cal S$, a contradiction.
So there is at least one edge between $T$ and $T'$ other than $x_1x'_1$.

Consider the case when the edges between $T$ and $T'$ are all incident with $x'_1$.
From above we know that $x'_1$ is adjacent to one of $w, y, z$.
By Lemma \ref{s4.5}, $x'_1$ is completely adjacent to $T$ by positive edges, except
possibly $x'_1w$ (which may or may not be an edge). If $x'_1w$ is a negative edge, 
then $T-\{x_1, x_2, \dots, x_{k-1}\}$ and $T'$ induce a graph in $\cal S$, 
a contradiction. So $x'_1$ is completely adjacent to $T$ by positive edges.  

Since $S$ minimally separates $\widehat{H_1}$ and $\widehat{H_2}$, $x_1$ has 
a neighbour in $\widehat{H_2}$. 
Let $u$ be a neighbour of $x_1$ in $\widehat{H_2}$. Suppose that $u$ is adjacent to
a vertex in $T'$. Then we know from above that $u$ must be adjacent to one of 
$w', y', z'$. By Lemma \ref{s4.5}, $u$ is completely adjacent to $T'$ by positive 
edges except possibly $uw'$. If $uw'$ is a negative edge, then 
$T$, $u$, and $T'-\{x'_1, x'_2, \dots, x'_{\ell-1}\}$ form an induced a graph
in $\cal S$. On the other hand if $u$ is completely adjacent to $T'$ by positive
edges, then $T$, $u$, and $T'$ form an induced graph in $\cal J$. 
Hence no neighbour of $x_1$ in $\widehat{H_2}$ is adjacent to any vertex in $T'$. 
Let $u_1u_2 \dots u_t$ be a shortest path in $\widehat{H_2}$ from a neighbour of 
$x_1$ to $T'-x'_1$. Then $t \geq 3$.  
If $u_{t-1}$ is not adjacent to any of $w', y', z'$, then 
$u_t = x'_j$ for some $j$. Choose such a $j$ to be the largest. Then
$T$, $u_1, u_2, \dots, u_{t-1}$, and $T'-\{x'_1, x'_2, \dots, x'_{j-1}\}$ induce 
a graph in $\cal S$, a contradiction. 
Thus $u_{t-1}$ is adjacent to at least one of $w', y', z'$.
Let $u'_t$ be the vertex in $T'$ closest to $x'_1$ which is a neighbour of $u_{t-1}$.
Then $x_1x'_1x'_2 \dots u'_tu_{t-1}u_{t-2} \dots u_1x_1$ is a cycle of length
$\geq 6$. This cycle cannot be induced as $\widehat{G}$ is $\cal C$-free.  
Since any chord in the cycle is incident with $x'_1$, $x'_1$ is completely adjacent 
to $u_1u_2 \dots u_{t-1}u'_t$. In particular, the distance between $x'_1$ and 
$u'_t$ in $T'$ is at most 2. Since $\widehat{G}$ does not contain any graph in $D$
as a induced subgraph, all chords in the cycle are positive.
Applying Lemma \ref{s4.5} to $T'$ and $u_{t-1}$ we conclude that $u_{t-1}$ is
completely adjacent to $T'$ by positive edges except possibly $u_{t-1}w'$. 
If $u_{t-1}w'$ is a negative edge then $T$, $u_1, u_2, \dots, u_{t-1}$, and
$T'-\{x'_1, x'_2, \dots, x'_{\ell-1}\}$ form an induced graph in $\cal S$. 
If $u_{t-1}$ is completely adjacent to $T'$ by positive edges then $T$,
$u_1, u_2, \dots, u_{t-1}$, and $T'$ form an induced graph in $\cal J$.  

Suppose that there is an edge between $x_1$ and $T'$ (other than $x_1x'_1$) 
and an edge between $x'_1$ and $T$  (other than $x_1x'_1$).
Again from above we know that $x_1$ is adjacent to one of $w', y', z'$ and that 
$x'_1$ is adjacent to at least one of $w, y, z$. By Lemma \ref{s4.5}, $x_1$ is 
completely adjacent to $T'$ by positive edges except possibly $x_1w'$, and 
$x'_1$ is completely adjacent to $T$ by positive edges except possibly $x'_1w$.
If $x_1$ is completely adjacent to $T'$ by positive edges and $x'_1$ 
is completely adjacent to $T$ by positive edges, or $T \in W_1$ and $T' \in W_1$ then 
$T$ and $T'$ form a graph in $\cal J$ as an induced subgraph. 
Suppose that $T \in W_1$ and $T' \notin W_1$. 
If $x_1w'$ is positive and $x'_1w$ is negative then $T$ and $T'$ induce a graph 
in $\cal J$.
If $x_1w'$ is negative then $T$ and 
$T' - \{x'_1, x'_2, \dots, x'_{\ell-1}\}$ induce a graph in $\cal S$.
A similar discussion applies to the case when $T \notin  W_1$ and $T' \in W_1$. 
Suppose that $T \notin W_1$ and $T' \notin W_1$. If $x_1w'$ is negative, then 
$T$ and $T'-\{x'_1, x'_2, \dots, x'_{\ell-1}\}$ induce a graph in $\cal S$.
Similarly, if $x'_1w$ is negative, then 
$T-\{x_1, x_2, \dots, x_{k-1}\}$ and $T'$ induce a graph in $\cal S$.

Suppose now that $x_1$ and $x'_1$ are in the same partite set of $\widehat{G}$.
Assume first that $x_1$ is adjacent to a vertex in $T'$ and $x'_1$ is adjacent to 
a vertex in $T$. We know from above that $x_1$ is adjacent to one of $w', y', z'$
and $x'_1$ is adjacent to one of $w, y, z$. By Lemma \ref{s4.5}, $x_1$ is 
completely adjacent to $T'$ by positive edges except possibly $x_1w'$, and
$x'_1$ is completely adjacent to $T$ by positive edges except possibly $x'_1w$.
If $x_1$ is completely adjacent to $T'$ by positive edges and $x'_1$ is completely 
adjacent to $T$ by positive edges then $T$ and $T'$ induce a graph in $\cal J$. 
If $x_1w'$ is a negative edge then $T$ and $T' - \{x'_1, x'_2, \dots, x'_{\ell-1}\}$
induce a graph in $\cal S$. Similarly, if $x'_1w$ is negative then 
$T-\{x_1, x_2, \dots, x_{k-1}\}$ and $T'$ induce a graph in $\cal S$. 

Assume next that $x_1$ is adjacent to a vertex in $T'$ but $x'_1$ is adjacent to
no vertex in $T$. As above, we know that $x_1$ is completely adjacent to $T'$ 
by positive edges except possibly $x_1w'$. If $x_1w'$ is a negative edge then $T$ and 
$T'-\{x'_1, x'_2, \dots, x'_{\ell-1}\}$ induce a graph in $\cal S$. So
$x_1$ completely adjacent to $T'$ by positive edges.

Let $u$ be a neighbour of $x'_1$ in $\widehat{H_1}$. Suppose that $u$ is adjacent to
a vertex in $T - x_1$. A similar proof as above shows that $u$ must be adjacent to
one of $w, y, z$. Observe that $T$ and $x'_1, x'_2$ form an induced tadpole with 
$x'_1$ being the end. Applying Lemma \ref{s4.5} to this tadpole and $u$ we conclude
that $u$ is completely adjacent to this tadpole by positive edges except possibly
$uw$ and $ux'_1$. If $uw$ is a negative edge then $T-\{x_1, x_2, \dots, x_{k-1}\}$,
$u$, and $T'$ form an induced graph in $\cal S$. So $u$ is completely adjacent to
$T$ by positive edges. Hence $T$, $u$, and $T'$ form an induced graph in 
$\cal J$. Hence no neighbour of $x'_1$ in $\widehat{H_1}$ is adjacent to any vertex
in $T - x_1$. 

Let $u_1, u_2, \dots, u_t$ be a shortest path in $\widehat{H_1}$ from a neighbour 
of $x'_1$ in $\widehat{H_1}$ to $T-x_1$. Then $t \geq 3$. 
A similar proof as above shows that $x_1$ is completely adjacent to 
$u_1u_2\dots u_{t-1}u'_t$ where $u'_t$ is the neighbour of $u_{t-1}$ in $T$ 
closest to $x_1$ by positive edges. Moreover, $u_{t-1}$ is completely adjacent to
$T$ by positive edges except possibly $u_{t-1}w$ and $u_{t-1}x_1$. 
Since $\widehat{G}$ contains no induced graph in $D$, $u_{t-1}x_1$ is a positive 
edge. If $u_{t-1}w$ is a negative edge, then $T-\{x_1,x_2,\dots,x_{k-1}\}$, 
$u_1, u_2, \dots, u_{t-1}$, $T'$ form an induced graph in $\cal S$. 
So $u_{t-1}$ is completely adjacent to $T$ by positive edges. 
Then $T$, $u_1, u_2, \dots, u_{t-1}$, and $T'$ form an induced graph in 
$\cal J$. 
A similar discussion applies when $x_1$ is adjacent to no vertex in $T'$ but $x'_1$ is adjacent to
a vertex in $T$.

It remains to consider the case when $x_1$ is adjacent to no vertex in $T'$ and 
$x'_1$ is adjacent to no vertex in $T$. By Lemma \ref{min-ss}, $x_1$ and $x'_1$ 
have a common neighbour $u$ in $\widehat{H_1}$ and a common neighbour $u'$ in 
$\widehat{H_2}$. If $u$ is not adjacent to any vertex of $T - x_1$ then 
$T$, $T'$ and $u$ induce a graph in $\cal S$. Similarly, if $u'$ is not adjacent to 
any vertex of $T' - x'_1$ then $T$, $T'$ and $u'$ induce a graph in 
$\cal S$. So $u$ is adjacent to some vertex in $T-x_1$ and $u'$ is adjacent to some
vertex in $T'-x'_1$. We know from above that $u$ must be adjacent to one of 
$w, y, z$. By Lemma \ref{s4.5} $u$ is completely adjacent to $T$ by positive edges
except possibly $uw$ and $ux_1$. If $uw$ is a negative edge then 
$T-\{x_1, x_2, \dots, x_{k-1}\}$, $u$ and $T'$ induce a graph in $\cal S$. 
If $ux_1$ is a negative edge then $\widehat{G}$ contains a graph in $D$ induced by
$x_1, u', x'_1, u, x_3, x_2$ when $k \geq 3$, or by 
$x_1, u', x'_1, u, y, x_2$ when $k = 2$, or by
$x_1, u', x'_1, u, w, y$ when $k = 1$. 
Hence $u$ is completely adjacent to $T$ by positive edges. 
A similar argument shows that $u'$ is completely adjacent to $T'$ by positive edges.
Therefore $T$ and $T'$ together with $u$ and $u'$ induce a graph in $\cal J$.
This completes the proof.
\qed

\begin{cor} \label{maincor}
A signed bigraph $\widehat{G}$ is not chordal if and only if it contains a 
non-trivial induced subgraph which has no signed simplicial edge.
\end{cor}
\pf Suppose that $\widehat{G}$ is not chordal. Then by Theorem \ref{main} 
$\widehat{G}$ contain a graph in $\cal F$ as an induced subgraph which is 
non-trivial and has no signed simplcial edge. On the other hand, suppose that 
$\widehat{G}$ is chordal. If $\widehat{H}$ is a non-trivial induced subgraph of
$\widehat{G}$, then the first edge in a signed simplicial edge-without-vertex 
ordering of $\widehat{G}$ that is contained in $\widehat{H}$ is a signed simplicial 
edge in $\widehat{H}$. 
\qed

\end{document}